\documentclass[12pt]{amsart}

\usepackage{amsmath,amsfonts,amsthm}
\usepackage{amscd}
\usepackage[psamsfonts]{amssymb}
\usepackage{graphics,epsfig,bbm,amssymb,amsmath}

\newtheorem{thm}{Theorem}[section]
\newtheorem{propos}{Proposition}[section]
\newtheorem{lem}{Lemma}[section]
\newtheorem{corrol}{Corollary}[section]

\newcounter{remarques}[section]
\renewcommand{\theremarques}{\arabic{section}.\arabic{remarques}}
\def\Rm{\refstepcounter{remarques}{\bf Remark \theremarques} }
\def\ali{\hfill\break}

\def\demi{{1\over 2}}

\def\Z{{{\Bbb Z}}}
\def\P{{{\Bbb P}}}
\def\N{{{\Bbb N}}}
\def\R{{{\Bbb R}}}

\def\aaa{{{\mathcal A}}}

\def\iii{{{\mathcal I}}}

\def\E{{{\Bbb E}}}

\def\indic{{{\bf 1}}}

\def\equilaw{{\stackrel{\hbox{law}}{=}}}
\def\zzz{{\mathcal Z}}
\def\Qneg{{Q^{\le 0}}}
\def\tQpos{{\tilde Q^{>0}}}
\def\intH{{\lfloor H\rfloor}}
\def\Qi{{Q^{\mathcal I}}}
\def\Ih{{\mathcal I_h^{(t)}}}
\def\tIh{{\tilde{\mathcal I}_h^{(t)}}}
\def\third{{1\over 3}}
\def\thQpos{{\hat Q^{>0}}}
\def\sss{{\mathcal S}}

\def\r{{\mathbb R}}
\def\e{{\mathbb E}}

\def\z{{\mathbb Z}}

\def\ee{\mathrm{e}}
\def\d{\, \mathrm{d}}

\def\ds{\displaystyle}
\def\eps{\varepsilon}
\def\R{{\mathbb R}}
\def\Z{{\mathbb Z}}

\def\N{{\mathbb N}}
\def\E{{\mathbb E}}

\def\t{\tilde}

\author[N. Enriquez]{Nathana\"el ENRIQUEZ}
\address{Laboratoire de Probabilit\'es et Mod\`eles Al\'eatoires, CNRS UMR 7599, Universit\'e Paris 6, 4
place Jussieu, 75252 Paris Cedex 05}
\email{enriquez@ccr.jussieu.fr}

\author[C. Sabot]{Christophe Sabot}
\address{Institut Camille Jordan, CNRS UMR 5208, Universit\'e de
Lyon 1, Universit\'e de Lyon, 43, Boulevard du 11 novembre 1918,
69622 Villeurbanne Cedex} \email{sabot@math.univ-lyon1.fr}

\author[O. Zindy]{Olivier ZINDY}
\address{Laboratoire de Probabilit\'es et Mod\`eles Al\'eatoires, CNRS UMR 7599, Universit\'e Paris 6, 4
place Jussieu, 75252 Paris Cedex 05} \email{zindy@ccr.jussieu.fr}


\title[Kesten's renewal constant]{A probabilistic representation of constants in Kesten's renewal theorem}

\begin{document}
\maketitle

\bigskip

{\footnotesize \noindent{\slshape\bfseries Abstract.} The aims of
this paper are twofold. Firstly, we derive a probabilistic
representation for the constant which appears in the one-dimensional
case of Kesten's renewal theorem. Secondly, we estimate the tail of
a related random variable which plays an essential role in the
description of  the stable limit law of one-dimensional transient
sub-ballistic random walks in random environment.}

\bigskip
\bigskip

\section{Introduction}
\label{sec:intro}

In 1973, Kesten published a famous paper \cite{kesten73} about the
tail estimates of renewal series of the form $\sum_{i\ge1} A_1\dots
A_{i-1} B_i,$ where $(A_i)_{i \ge 0}$ is a sequence of non-negative
i.i.d. $d\times d$ random matrices and $(B_i)_{i \ge 1}$ is a
sequence of i.i.d. random vectors of $\r^d.$ His result states that
the tail of the projection of this random vector on every direction
is equivalent to $C t^{-\kappa},$ when $t$ tends to infinity, where
$C$ and $\kappa$ are positive constants. The constant $\kappa$ is
defined as the solution of the equation $k(s)=1,$ with $k(s):=
\lim_{n \to \infty} \e(\parallel A_1 \dots A_n \parallel^s)^{1/n}.$
The proof of his result in the one-dimensional case, even if it is
much easier than in dimension $d \ge 2,$ is already rather
complicated.

Even though we are concerned by the one-dimensional case in this paper,
let us mention that a significant generalization of Kesten's result,
in the multi-dimensional case, was recently achieved by de Saporta,
Guivarc'h and Le Page \cite{desaporta-guivarch-lepage}, who relaxed
the assumption of positivity on $A_i.$

In 1991, Goldie \cite{goldie} relaxed, in dimension $d=1,$ the
assumption of positivity on the $A_i$ and simplified Kesten's proof.
Furthermore, he obtained a formula for the implicit constant $C$ in
the special case where $A_i$ is non-negative and $\kappa$ is an
integer.

In 1991, Chamayou and Letac \cite{chamayou-letac} observed that, in
dimension $d=1,$ if $A_i$ has the same law as $(1-X_i)/X_i,$ with
$X_i$ following a Beta distribution on $(0,1),$ then the law of the
series itself is computable so that the constant $C$ is explicit in
this special case also. The following question was then asked. How
does one effectively compute the constant $C$?

In our framework, we consider the case $d=1$ and we make the
following assumptions: $\rho_i=A_i$ is a sequence of i.i.d. positive
random variables, $B_i=1$ and there exists $\kappa>0$ such that
$\e(\rho_1^\kappa)=1.$ Moreover, we assume a weak integrability
condition and that the law of $\log \rho_i,$ which has a negative
expectation by the previous assumptions, is non-arithmetic. In this
context we are interested in the random series
$$
R=1+\sum_{k\ge 1} \rho_1\cdots \rho_k.
$$
The previous assumptions ensure that the tail of the renewal series $R$
is equivalent to $C_K t^{-\kappa},$ when $t$ tends to infinity. We
are now aiming at finding a probabilistic representation of the
constant $C_K$.

Besides, this work is motivated by the study of one-dimensional
random walks in random environment. In \cite{kesten-kozlov-spitzer},
Kesten, Kozlov and Spitzer proved, using the tail estimate derived
in \cite{kesten73}, that when the RWRE is transient with null
asymptotic speed, then the behavior depends on an index $\kappa\le
1$: the RWRE $X_n$ normalized by $n^{1/\kappa}$ converges in law to
$C_\kappa \left( {1\over \sss_\kappa}\right)^\kappa$ where
$\sss_\kappa$ is a positive stable random variable with index
$\kappa.$ The computation of the explicit value of $C_\kappa$ was
left open. In \cite{enriquez-sabot-zindy}, the authors derive an
explicit expression, either in terms of the Kesten's constant $C_K$
when it is explicit, or in terms of the expectation of a random
series when $C_K$ is not explicit. To this end, we need to obtain a
tail estimate for a random variable $Z$, closely related to the random
series $R$, and to relate it to Kesten's constant. This is the other
aim of this paper.

The strategy of our proof is based on a coupling argument in the
(cf \cite{durrett}, 4.3). We first interpret $\rho_1\dots \rho_n$ as
the exponential of a random walk $(V_n \, ,\, n \ge 0),$ which is
negatively drifted, since $\e(\log \rho_1)<0.$ We have now to deal
with the series $R:=\sum_{n \ge 0} \ee^{V_n}.$ One can write
$$
R=\ee^{S} \sum_{n \ge 0} \ee^{V_n-S},
$$
where $S$ is the maximum of $(V_n \, ,\, n \ge 0).$ The heuristic is
that $S$ and  $\sum_{n \ge 0} \ee^{V_n-S}$ are asymptotically
independent. The coupling argument is used to derive this asymptotic
independence. But, in order to implement this strategy, several
difficulties have to be overcome: we first need to condition $S$ to
be large. Moreover, we have to couple conditioned processes: this
requires us to describe precisely the part of the process $(V_0, \ldots
,V_{T_S})$, where $T_S$ is the first hitting time of the level $S$.

To end this section, let us finally discuss our results and
strategy. Let us first remind that Kesten and Goldie's proof were
based on a clever use of the renewal theorem but strongly relied on
the renewal structure of the series, and also did not lead to satisfying
representations of the constant involved in its tail function.
Later, Siegmund \cite{siegmund} presented an interesting scheme of
proof,  
inspired by a  work on change-point analysis of Pollak and Yakir \cite{yakir-pollak}. 
He was able to derive formally a representation of the constant, which enables simulation
of the constant by Monte Carlo. 

We would like to emphasize the flexibility of our proof that allows
to study conditioned variables which do not necessarily satisfy a
renewal scheme like the variable $Z$ mentioned above, which plays a key role in
the analysis of RWRE. This flexibility could hopefully make also
possible some generalizations to the $d$-dimensional case. As explained above,
the strength of this method is indeed to prove an asymptotic
independence between  two different parts of the underlying random
walk of step $\log (\rho_n)$, when its maximum is large, namely : the
maximum of the random walk  and the part of 
trajectory in the neighbourhood  of the absolute maximum.
As a consequence, the tail constant of $R$ is expressed as the product of the tail constant of the absolute maximum of the random walk times the expectation of a functional of some random walk which comes from the part of the trajectory near its maximum. 
One of the central interests of this representation is that it is well suited for Monte-Carlo simulation.
Compared to Siegmund's formula, our formula is exact and not asymptotic (formula (3.6) of
Siegmund \cite{siegmund} must be understood as a limit when $j$ tends to infinity). Our asymptotic independence argument  is reminiscent of the argument of Siegmund which remained at a heuristic level, and we want to emphasize that this asymptotic independence is the difficult part of our proof.

On the other hand, let us notice that the analytic expressions found
by Goldie when $\kappa$ is an integer, and Chamayou and Letac in the
case of Beta variables are strongly based on the renewal scheme. It
is therefore not surprising that the representation found by our
method do not recover these results. However, their identification a
posteriori leads to explicit formulas for the constants arising in
the limit theorems for RWRE in some very interesting special cases
see \cite{enriquez-sabot-zindy}.

\section{ Notation and statement of the results}
\label{sec:main}
 Let $(\rho_i)_{i\in \Z}$ be a sequence of i.i.d.
positive random variables with law $Q=\mu^{\otimes \Z}$. With the
sequence $(\rho_i)_{i\in \Z}$ we associate the potential
$(V_k)_{k\in \Z}$ defined by
$$
V_n :=\left\{\begin{array}{lll} \sum_{k=1}^n \log \rho_k & {\rm if}
\ n \ge 1,
\\
 0 &  {\rm if} \ n=0,
\\
-\sum_{k=n+1}^0 \log \rho_k &{\rm if} \ n\le -1.
\end{array}
\right.
$$

Let $\rho$ have law $\mu.$ Suppose now that the law $\mu$ is such
that there is $\kappa>0$ satisfying
\begin{eqnarray}
\label{kappa} \E^\mu(\rho^\kappa)=1 \qquad {\rm and} \qquad \E^\mu(
\rho^{\kappa} \log^+ \rho)< \infty.
\end{eqnarray}

\noindent Moreover, we assume that the distribution  of $\log \rho$
is non-lattice. Then the law $\mu$ is such that $\log \rho$ satisfies
\begin{eqnarray}\label{transience}
\E^\mu(\log \rho)<0,
\end{eqnarray}
which implies that, $Q$-almost surely,
$$
\lim_{n\to\infty} {V_n\over n}=\int \log\rho \d\mu <0.
$$

We set
$$
S:=\max\{ V_k, \; k\ge 0\},
$$
and
$$
H:=\max\{V_k, 0\le k\le T_{\R_-}\},
$$
where $T_{\R_-}$ is the first positive hitting time of $\R_-$:
$$
T_{\R_-}:=\inf\{k>0, \; V_k\le 0\}.
$$
The random variable $S$ is the absolute maximum of the path
$(V_k)_{k\ge 0}$ while $H$ is the maximum of the first positive
excursion. We also set
$$
T_S:=\inf\{k\ge 0, \; V_k=S\}, \;\;\; T_H:=\inf\{k\ge 0, \; V_k=H\}.
$$
We clearly have, $Q$-almost surely,
$$
H\le S<\infty, \;\;\; T_H\le T_S<\infty.
$$

The following tail estimate for $S$ is a classical consequence of
renewal theory, see \cite{feller},
\begin{eqnarray}\label{Feller}
\P^Q(\ee^S\ge t)\sim C_F t^{-\kappa},
\end{eqnarray}
when $t\to \infty$, where
$$C_F={1-\e^Q(\ee^{\kappa V(T_{\R_-})})\over \kappa
\e^\mu(\rho^\kappa\log\rho)\e^Q(T_{\R_-})}.
$$

\noindent The tail estimate of $H$ is derived by Iglehart, in
\cite{igle},
\begin{eqnarray}\label{Igle}
\P^Q(\ee^H\ge t)\sim C_I t^{-\kappa},
\end{eqnarray}
when $t\to \infty$, where
\begin{eqnarray*}\label{CI}
C_I={(1-\e^Q(\ee^{\kappa V(T_{\R_-})}))^2\over \kappa
\e^\mu(\rho^\kappa\log\rho)\e^Q(T_{\R_-})}=(1-\e^Q(\ee^{\kappa
V(T_{\R_-})}))C_F.
\end{eqnarray*}

Consider now the random variable
$$
R:=\sum_{n=0}^\infty \ee^{V_n}.
$$
This random variable clearly satisfies the following random affine
equation
$$
R\stackrel{\hbox{law}}{=} 1+\rho R,
$$
where $\rho$ is a random variable with law $\mu$ independent of $R$.
In \cite{kesten73}, Kesten proved (actually his result was more
general and concerned by the multidimensional version of this one)
that there exists a positive constant $C_K$ such that
\begin{eqnarray}
\label{equivkesten} \P^Q(R\ge t)\sim C_K t^{-\kappa},
\end{eqnarray}
when $t\to \infty$. The constant $C_K$ has been made explicit in
some particular cases: for $\kappa$ integer by Goldie, see
\cite{goldie}, and when $\rho \stackrel{\hbox{law}}{=} {W\over 1-W}$
where $W$ is a beta variable, by Chamayou and Letac
\cite{chamayou-letac}. One aim of this paper is to derive an
expression of this constant in terms of the expectation of a functional
of the random walk $V$ which is more standard than $R.$

We need now to introduce some Girsanov transform of $Q$. Thanks to
(\ref{kappa}) we can define the law
$$
\tilde\mu=\rho^\kappa \mu,
$$
and the law $\tilde Q=\tilde\mu^{\otimes \Z}$ which is the law of a
sequence of i.i.d. random variables with law $\tilde\mu$. The
definition of $\kappa$ implies that
$$
\int \log\rho \, \tilde\mu(\d\rho) >0,
$$
and thus that, $\tilde Q$-almost surely,
$$
\lim_{n\to\infty} {V_n\over n}=\int \log \rho \, \d\tilde \mu >0.
$$
 Moreover, $\tilde Q$ is a Girsanov transform of
$Q$, i.e. we have for all $n$
$$
\E^Q\left(\phi(V_0, \ldots , V_n)\right)=\E^{\tilde Q} \left(
e^{-\kappa V_n} \phi(V_0, \ldots ,V_n)\right),
$$
for any bounded test function $\phi$.
 Let us now introduce the random variable $M$ defined by
 \begin{eqnarray} \label{M}
 M=\sum_{i<0} \ee^{-V_i}+ \sum_{j\ge 0} \ee^{-V_j},
 \end{eqnarray}
 where $(V_i)_{i<0}$ is distributed under $Q(\cdot | V_i\ge 0,
 \forall i<0)$ and independent of $(V_j)_{j\ge 0}$ which is
 distributed under $\tilde Q(\cdot| V_j>0, \forall j> 0)$.
\medskip

\begin{thm}\label{Kesten}
i) We have the following tail estimate
$$
\P^Q(R\ge t)\sim C_K t^{-\kappa},
$$
when $t \to \infty,$ where
$$
C_K = C_F \E(M^\kappa).
$$

ii) We have
$$
\P^Q(R\ge t\, ; \,  H=S)\sim C_{KI} t^{-\kappa},
$$
when $t \to \infty,$ where
$$
C_{KI} := C_I \E(M^\kappa).
$$
\end{thm}

\medskip
\noindent\Rm: The conditioning $H=S$ means that the path $(V_k)_{k
\ge 0}$ never goes above the height of its first excursion. \ali

 In \cite{enriquez-sabot-zindy}, we need a tail estimate on a random variable of the type of $R$ but
 with an extra term. Let us introduce the event
 \begin{eqnarray} \label{iii}
  \iii:=\{H=S\}\cap \{V_k\ge 0 \, , \, \forall k\le 0\},
 \end{eqnarray}
and the random variable
$$
 Z:=\ee^S M_1 M_2,
$$
 where
 $$
 M_1:= \sum_{k=-\infty}^{T_S} \ee^{-V_k},
 $$
 $$
 M_2:=\sum_{k=0}^\infty \ee^{V_k-S}.
 $$
 \medskip

 \begin{thm}\label{doubleM}
 We have the following tail estimate
 $$
 \P^Q( Z\ge t|\iii)\sim {1\over \P^Q(H=S)}C_U t^{-\kappa},
 $$
when $t \to \infty,$ where
 $$
 C_U=C_I\E(M^\kappa)^2={C_{I}\over C_F} (C_K)^2.
 $$
 \end{thm}
\medskip

\noindent \Rm:  The conditioning event $\iii$ gives a nice symmetry property, which is useful to return the path, cf Subsection \ref{subsec:return}.

Let us now discuss the case where the $B_i$'s are not necessarily
equal to $1.$ Let $(B_i)_{i \ge 0}$ be a sequence of positive i.i.d.
random variables, which is independent of the sequence $(\rho_i)_{i
\ge 0},$ and denote by $R^B$ the random series $R^B:=B_0+\sum_{k\ge
1} B_k \rho_1 \cdots \rho_k .$ The result of Theorem \ref{Kesten},
$i),$ is then generalized into the following result.
\medskip
 \begin{thm}\label{Kestengene}
 If there exists $\varepsilon>0$ such that $\E(|B_1|^{\kappa+ \varepsilon})<
 \infty,$ then
$$
\P^Q(R^B\ge t)\sim C_{KB} t^{-\kappa},
$$
when $t \to \infty,$ where
$$
C_{KB} = C_F \E((M^B)^\kappa)
$$
and where $M^B$ is defined by
$$
M^B=\sum_{k<0} \ee^{-V_k} \tilde{B}_k+ \sum_{k\ge 0}
\ee^{-V_k}\tilde{B}_k,
$$
with $(V_k)_{k<0}$ distributed under $Q(\cdot | V_k\ge 0,
 \forall k<0)$ and independent of $(V_k)_{k\ge 0}$ which is
distributed under $\tilde Q(\cdot| V_k>0, \forall k> 0)$ while
$(\tilde{B}_k)_{k \in \z}$ is a sequence of i.i.d. random variables
having the same distribution as $B_1$ and independent of
$(V_k)_{k\in \z}.$
\end{thm}
\medskip

 \ali\ali
{\bf  Sketch of the proof and organization of the paper}
 \ali

The intuition behind Theorem \ref{Kesten} and Theorem \ref{doubleM}
is the following. Let us first consider $\P^Q(R\ge t | H=S)$. The
law $Q(\cdot| \iii)$ has a symmetry property which implies that the
variable $R=M_2 \ee^H$ has the same distribution as $M_1 \ee^H$ (cf
Subsection \ref{subsec:return}). Then, the proof of Theorem \ref{Kesten}
is based on the following  arguments.

Firstly, we prove that the variables $M_1$ and $\ee^H$ are
asymptotically independent. To this end, we use a delicate coupling
argument which works only when $H$ is conditioned to be large.
Therefore, we need to restrict ourselves to large values of $H$. To
this end, we need to control the value of $R$ conditioned by $H$;
this is done in Section \ref{sec:prelimin}. Then, a second difficulty is that we have
to couple conditioned processes (namely, the process $(V_k)$
conditioned to have a first high excursion). We overcome this
difficulty by using an explicit description of the law of the path
$(V_0, \ldots, V_{T_H})$. Namely, the  path $(V_0, \ldots ,V_{T_H})$
behaves like $V$ under $\tilde Q(\cdot| V_k>0, \forall k>0)$ stopped
at some random time.

Secondly, we observe that the distribution of $M_1$  is close to the
distribution of $M$ as a consequence of the above description of the
law of  $(V_0, \ldots ,V_{T_H})$.

From these two facts, we deduce that $\P^Q(R\ge t\;|\;
\iii)\simeq\P^Q(M \ee^H\ge t\;|\; \iii)$, where $M$ and $H$ are
roughly independent. Using the tail estimate for $H$ we get the part
ii) of Theorem \ref{Kesten}. For Theorem \ref{doubleM}, we proceed
similarly: the variable $Z$ can be written $M_1R$ and,
 for large $H$, the variables $M_1$ and $R$ are asymptotically
 independent and the law of $M_1$ is close to the law
of $M$. Then the estimate on the tail of $R$ allows us to conclude the proof.

Let us now describe the organization of the proofs. Section \ref{sec:preliminaries} contains preliminary results, whose proofs are postponed to the Appendix (see Section \ref{sec:appendix}). In Subsection
\ref{subsec:moment}, we prove that $M$ has finite moments of all orders
and we estimate the rest of the series $M$. Subsection \ref{subsec:return}
contains some preliminary properties of the law $Q(\cdot|\iii)$, and
Subsection \ref{subsec:artificial} presents a representation of the law of
the process $(V_0, \ldots ,V_S)$ in terms of the law $\tilde Q$.
Section \ref{sec:prelimin} contains crucial estimates which will
allow us to restrict ourselves to large values of $H$. In Section
\ref{sec:coupling}, we detail the coupling arguments which roughly
give the asymptotic independence of $M_1 $ and $\ee^H M_2$. Finally,
in Section \ref{sec:proof} we assemble the arguments of the previous
sections to prove Theorem \ref{doubleM} and Theorem \ref{Kesten}. In the Appendix (see
Section \ref{sec:appendix}), we give the proof of the claims of Section \ref{sec:preliminaries} and present a Tauberian version of the
tail estimates, which is the version we ultimately use in
\cite{enriquez-sabot-zindy}.

Let us finally explain the convention we use concerning constants. We denote by $c$ a positive constant with value changing from place to place, which only depends on $\kappa$ and the distribution of $\rho.$  The dependence on additional parameters otherwise appears in the notation.

\section{Preliminaries}
\label{sec:preliminaries}
In this section, we give preliminary results, whose proof are postponed to the Appendix (see Section \ref{sec:appendix}).

\subsection{Moments of $M$}
\label{subsec:moment}
Here is a series of three lemmas about the moments of the exponential
functional of the random walk $M.$ In this section, we denote by $\{V \ge
-L\}$ the event $\{V_k \ge -L, \, k \ge 0\}.$
\medskip
\begin{lem}\label{l:moment1}
 There exists $c>0$ such that, for all $L\geq0$,
$$\e^{\t Q}\Big(\sum_{k\geq0} \ee^{-V_k}\,|\, V\geq -L\Big)\le c \, \ee^L.$$
\end{lem}
\medskip
\begin{lem}
\label{l:moment2} Under $\t Q^{\ge0}:=\t Q(\cdot \, | \, V_k\geq0 \,
, \, \forall k\geq0),$ all the moments of $\sum_{k\geq 0}
\ee^{-V_k}$ are finite.
\end{lem}
\medskip

We will need further a finer result than Lemma \ref{l:moment1} as follows.
\medskip
\begin{lem}
\label{l:moment3}  For any $\kappa>0,$ there exists $c=c(\kappa)>0$ such that, for all $L>0$ and
for all $\eps'>0$, we have
\begin{itemize}
  \item if $\kappa <1,$
  $$\e^{\t Q}\Big(\sum_{i\geq0} \ee^{-V_i}\, |\, V\geq -L\Big)\leq c \, \ee^{(1-\kappa+\eps')L},$$
  \item if $\kappa \ge 1,$
  $$\e^{\t Q}\Big(\sum_{i\geq0} \ee^{-V_i}\, |\, V\geq -L\Big)\leq c \, \ee^{\eps'L}.$$
\end{itemize}
\end{lem}
\medskip

\noindent\Rm: Analogous results as in Lemma \ref{l:moment1}, Lemma
\ref{l:moment2} and Lemma \ref{l:moment3} apply for $\sum_{k\geq 0}
\ee^{V_k}$ under $Q$ and conditionally on the event $\{ V_k\leq L \,
, \, \forall k\ge0\}$. \ali

\subsection{A time reversal}
\label{subsec:return} Let us denote by $Q^\iii$ the conditional law
$Q^\iii(\cdot):=Q(\cdot|\iii),$ where $\iii$ is defined in
(\ref{iii}). The law $Q^\iii$ has the following symmetry property.
\medskip

\begin{lem}
\label{l:return}
 Under $Q^\iii$ we have the following equality in law
$$
(V_k)_{k\in \Z}\stackrel{\hbox{law}}{=}(V_{T_H}-V_{T_H-k})_{k\in
\Z}.
$$
\end{lem}
\medskip

This implies that under $Q^\iii$, $R$ has the law of $\ee^H M_1$.
This last formula will be useful since the asymptotic independence
of $\ee^H$ and $M_1$, in the limit of large $H$, is more visible
than the asymptotic independence of $H$ and $M_2$ and will be easier
to prove.

\subsection{The two faces of the mountain}
\label{subsec:artificial}
 It will be convenient to introduce the
following notation: we denote by $Q^{\le 0}$ the conditional law
$$
\Qneg(\cdot)=Q(\cdot | V_k\le 0, \; \forall k\ge 0),
$$
and by $\tQpos$ the conditional law
$$
\tQpos(\cdot)=\tilde Q(\cdot | V_k>0, \;\forall k>0).
$$

It will be useful to describe the law of the part of the path $(V_0,
\ldots, V_{T_S})$. Let us introduce some notations. If $(Y_k)_{k\ge
0}$ is a random process under the law $\tilde Q$, then $Y_k\to
+\infty$ a.s. and we can define its strictly increasing ladder times
$(e_k)_{k\ge 0}$ by: $e_0:=0$, and
$$
e_{k+1}:=\inf\{n>e_k, \; Y_n>Y_{e_k}\}.
$$
We define a random variable $((Y_k)_{k\ge 0}, \Theta)$ with values
in $\R^\N\times \N$ as follows: the random process $(Y_k)_{k\ge 0}$
has a law with density with respect to $\tilde Q$ given by
$$
{1\over \zzz} \Big(\sum_{k=0}^\infty \ee^{-\kappa Y_{e_k}}\Big)
\tilde Q,
$$
where $\zzz$ is the normalizing constant given by
$$
\zzz={1\over 1-\E^{\tilde Q}(\ee^{-\kappa Y_{e_1}})}.
$$
Then, conditionally on $(Y_k)_{k\ge 0}$, $\Theta$ takes one of the
value of the strictly ladder times with probability
$$
\P(\Theta=e_p\;|\; \sigma((Y_k)_{k\ge 0}))={\ee^{-\kappa
Y_{e_p}}\over\sum_{k=0}^\infty \ee^{-\kappa Y_{e_k}}}.
$$
We denote by $\hat Q$ the law of $((Y_k)_{k\ge 0}, \Theta)$.
Otherwise stated, it means that, for all test functions $\phi,$
$$
\E^{\hat Q} (\phi(\Theta, (Y_n)_{n\ge 0}))={1\over \zzz}\E^{\tilde
Q}\Big( \sum_{k=0}^\infty \ee^{-\kappa Y_{e_k}} \phi(e_k, (Y_n)_{n
\ge 0})\Big).
$$
\medskip

\begin{lem}\label{artificial}
The processes $(V_0, \ldots, V_{T_S})$ and
$(V_{T_{S+k}}-V_{T_S})_{k\ge 0}$ are independent and have the
following laws: $(V_{T_{S+k}}-V_{T_S})_{k\ge 0}$ has the law $Q^{\le
0}$ and
$$
(V_0, \ldots ,V_{T_S})\equilaw (Y_0, \ldots ,Y_\Theta),
$$
where $((Y_k)_{k\ge 0},\Theta)$ has the law $\hat Q$.
\end{lem}
\medskip

Denote now by $\hat Q^{>0}$ the law
$$
\thQpos =\hat Q(\cdot \; | \; Y_k>0, \forall k>0).
$$
We will need the following result.
\medskip

\begin{lem}\label{artificial2}
There exists a positive constant $c>0$ such that, for all positive
test functions $\psi,$
$$
\E^{\Qi}(\psi(V_0, \ldots ,V_{T_H}))\le c \E^{\thQpos}( \psi(Y_0,
\ldots , Y_\Theta)).
$$
\end{lem}
\medskip

\section{A preliminary estimate}
\label{sec:prelimin} To derive the tail estimate of $R$ or $Z$ we
need to restrict to large values of $H$: this will be possible,
thanks to the following estimate.
\medskip

\begin{lem} \label{prelim}
For all $\eta>0$ there exists a positive constant $c_\eta$ such that
$$
\E^{Q^\iii}((M_1)^\eta\; |\; \lfloor H\rfloor)\le c_\eta , \;\;\;
Q^\iii \hbox{- a.s.,}
$$
where $\lfloor H\rfloor $ is the integer part of $H$.
\end{lem}

\medskip

\begin{proof} Since $(V_k)_{k\le 0}$  is independent of $H$
under $Q^\iii$, we have, for all $p\in \N,$
$$
\E^{Q^\iii}((M_1)^\eta\; |\; \lfloor H\rfloor=p)\le 2^\eta \left(
\E^{\Qneg}\bigg(\big(\sum_{k=0}^\infty \ee^{V_k}\big)^\eta \bigg)+
\E^{Q^\iii}\bigg(\big(\sum_{k=0}^{T_H} \ee^{-V_k}\big)^\eta \;|\;
\intH =p\bigg) \right).
$$
The first term on the right-hand side is finite for all $\eta>0$ as
proved in Subsection \ref{subsec:moment}. Consider now the last term.
Using Lemma \ref{artificial2}, we get
\begin{eqnarray*}
&&\E^{Q^\iii}\bigg(\big(\sum_{k=0}^{T_H} \ee^{-V_k}\big)^\eta \;|\;
\intH =p\bigg)
\\
&\le& {c\over \P^{Q^\iii}(\intH =p)}
\E^{\thQpos}\bigg(\big(\sum_{k=0}^{T_H} \ee^{-Y_k}\big)^\eta
\indic_{\intH =p}\bigg)
\\
&\le & {c \over \P^{Q^\iii}(\intH =p)} \E^{\tQpos}\bigg(\big(
\sum_{k=0}^\infty \ee^{-\kappa Y_{e_k}} \indic_{Y_{e_k}\in
[p,p+1[}\big) \big(\sum_{j=0}^{\infty} \ee^{-V_j}\big)^\eta \bigg).
\end{eqnarray*}
Now, using the Cauchy-Schwarz inequality in the last expression, we get
\begin{eqnarray*}
&&\E^{Q^\iii}\bigg(\big(\sum_{k=0}^{T_H} \ee^{-V_k}\big)^\eta \;|\;
\intH =p\bigg)
\\
&\le& {c\over \P^{Q^\iii}(\intH =p)} \E^{\tQpos}\bigg(\big(
\sum_{k=0}^\infty \ee^{-\kappa Y_{e_k}} \indic_{Y_{e_k}\in
[p,p+1[}\big)^2\bigg)^\demi \E^{\tQpos}\bigg(
\big(\sum_{k=0}^{\infty} \ee^{-V_k}\big)^{2\eta} \bigg)^\demi
\\
&\le& {c \, \ee^{-\kappa p}\over \P^{Q^\iii}(\intH =p)}
\E^{\tQpos}\bigg(\big( \sum_{k=0}^\infty \indic_{Y_{e_k}\in
[p,p+1[}\big)^2\bigg)^\demi \E^{\tQpos}\bigg(
\big(\sum_{k=0}^{\infty} \ee^{-V_k}\big)^{2\eta} \bigg)^\demi.
\end{eqnarray*}
But the last term is independent of $p$ and finite by Lemma
\ref{l:moment2}. On the other hand, since $\tilde{Q}(V_k>0 \, , \,
\forall k >0)>0$ and from the Markov property, we obtain
$$
\E^{\tQpos}\bigg(\big( \sum_{k=0}^\infty \indic_{Y_{e_k}\in
[p,p+1[}\big)^2\bigg) \le c \E^{\tilde{Q}}\bigg(\big(
\sum_{k=0}^\infty \indic_{Y_{e_k}\in [p,p+1[}\big)^2\bigg) \le c
\E^{\tilde{Q}}\bigg(\big( \sum_{k=0}^\infty \indic_{Y_{e_k}\in
[0,1[}\big)^2\bigg),
$$
which is finite since $(Y_k)_{k \ge 0}$ has a positive drift under
$\tilde Q$. Finally, using the tail estimate on $H,$ we know that
\begin{eqnarray}
\label{limiteHinteger} \lim_{p\to\infty} \ee^{\kappa p}
\P^{Q^\iii}(\intH =p) &=& \lim_{p\to\infty} \ee^{\kappa p}
\left(\P^{Q^\iii}(H \ge p)-\P^{Q^\iii}(H\ge p+1)\right)
\\
&=& C_I(1-\ee^{-\kappa}). \nonumber
\end{eqnarray}
Hence, $(\ee^{\kappa p} \P^{Q^\iii}(\intH =p))^{-1}$ is a bounded
sequence (we do not have to consider the cases where eventually
$\P(\intH=p)=0$ since it is a conditioning by an event of null
probability which can be omitted).
\end{proof}

\medskip

\begin{corrol}\label{corro61}
We have, $\Qi$-almost surely,
$$
\E^{\Qi}\left( Z\; |\; \intH \right) \le \ee c_2 \ee^{\intH}.
$$
\end{corrol}
\begin{proof} We have $Z=M_1 M_2 \ee^H$. Using the Cauchy-Schwarz inequality and
Lemma \ref{prelim} we get
\begin{eqnarray*}
\E^{\Qi}\left( Z\; |\; \intH \right)&\le& \ee^{\intH+1}
\left(\E^{\Qi}((M_1)^2\; |\; \intH)\E^{\Qi}((M_2)^2\; |\;
\intH)\right)^\demi
\\
&\le& \ee c_2 \ee^\intH,
\end{eqnarray*}
since $M_1$ and $M_2$ have the same law under $\Qi$.
\end{proof}

\medskip

\begin{corrol}\label{ht}
Let $h:\R_+\mapsto \R_+$ be a function such that
$$
\lim_{t\to \infty} t^{-1} \ee^{h(t)}=0.
$$
Then, we have
$$
\P^{Q^\iii}\left( R\ge t, \; H\le h(t)\right)=o(t^{-\kappa}),
$$
$$
\P^{Q^\iii}\left( Z\ge t, \; H\le h(t)\right)=o(t^{-\kappa}),
$$
when $t$ tends to infinity.
\end{corrol}

\medskip

\begin{proof}
Let us do the proof for $Z$. Let $\eta$ be a positive real such that
$$
\eta>\kappa.
$$
We have (all expectations are relative to the measure $Q^\iii$; so,
to simplify the reading, we remove the reference to $Q^\iii$ in the
following)
\begin{eqnarray*}
\P^{Q^\iii}\left( Z\ge t, \; H\le h(t)\right) &=& \E\left( \P\left
(Z\ge t, \; H\le h(t)\; | \; \intH \right)\right)
\\
&\le& \E\left( \indic_{\intH\le \lfloor h(t)\rfloor} \P\left (Z\ge t
\; | \; \intH \right)\right)
\\
&\le& \E\left( \indic_{\intH\le \lfloor h(t)\rfloor} \P\left (M_1
M_2 \ge t \ee^{-(\intH+1)} \; | \; \intH \right)\right)
\\
&\le&\ee^\eta \E\left( \indic_{\intH\le \lfloor h(t)\rfloor}
t^{-\eta} \ee^{\eta \intH}\E\left ((M_1 M_2)^\eta \; | \; \intH
\right)\right)
\\
&\le& \ee^\eta\E\left( \indic_{\intH\le \lfloor h(t)\rfloor}
t^{-\eta} \ee^{\eta \intH}\E\left ((M_1)^{2\eta} \; | \; \intH
\right)\right).
\end{eqnarray*}
In the last formula, we used the Cauchy-Schwarz inequality and the
symmetry property of $Q^\iii,$ see Lemma \ref{l:return}, to obtain
$$
\E((M_2)^{2\eta}\;|\; \intH)=\E((M_1)^{2\eta}\;|\; \intH).
$$
We can now use the estimate of Lemma \ref{prelim}, which gives
\begin{eqnarray*}
\P^{Q^\iii}\left( Z\ge t, \; H\le h(t)\right) &\le& \ee^\eta
c_{2\eta}t^{-\eta} \sum_{p=0}^{\lfloor h(t)\rfloor} \ee^{\eta p}
\P(\intH =p)
\\
&\le & c t^{-\eta} \sum_{p=0}^{\lfloor h(t)\rfloor}
\ee^{(\eta-\kappa) p}.
\end{eqnarray*}
In the last formula, we used the fact that
$\P(\intH=p)=O(\ee^{-\kappa p}),$ see (\ref{limiteHinteger}). Since
we chose $\eta>\kappa$ we can bound uniformly
$$
\P^{Q^\iii}\left( Z\ge t, \; H\le h(t)\right) \le c t^{-\eta }
\ee^{(\eta-\kappa)h(t)}= c t^{-\kappa}\Big({\ee^{h(t)}\over
t}\Big)^{\eta-\kappa}.
$$
This gives the result for $Z$. Since $R\le Z,$ we get the result for
$R$.
\end{proof}

\section{The coupling argument}
\label{sec:coupling} We set
$$
I(t):= \P^\Qi\left( \ee^H M_1 M_2\ge t\right),
$$
$$
J(t):= \P^\Qi\left( \ee^H M_2 \ge t\right),
$$
$$
K(t):= \P^\Qi\left( \ee^H \ge t\right).
$$
\noindent  From the estimate of Iglehart, see \cite{igle}, we know
that
$$
K(t)\sim {1\over \P^Q(H=S)} C_I t^{-\kappa},
$$
when $t \to \infty.$ Indeed, we have
$$
\P^\Qi\left( \ee^H \ge t\right)={1\over \P^Q(H=S)}(\P^Q(\ee^H\ge
t)-P^Q(\ee^H\ge t, \; S>H)).
$$
The second term is clearly of order $O(t^{-2\kappa})$, the first
term is estimated in \cite{igle}, cf (\ref{Igle}).

 We will prove the following key estimates.
\medskip
\begin{propos}\label{p:encadrement}
For all $\xi>0$ there exists a function $\epsilon_\xi(t)>0$ such
that $\lim_{t\to \infty} \epsilon_\xi(t)=0$
 and
$$
\ee^{-3\xi}\E\left( J(\ee^{3\xi} t
M^{-1})\right)(1-\epsilon_\xi(t))\le I(t)\le \ee^{3\xi}\E\left(
J(\ee^{-3\xi} t M^{-1})\right)(1+\epsilon_\xi(t)),
$$
$$
\ee^{-2\xi}\E\left( K(\ee^{2\xi} t
M^{-1})\right)(1-\epsilon_\xi(t))\le J(t)\le \ee^{2\xi}\E\left(
K(\ee^{-2\xi} t M^{-1})\right)(1+\epsilon_\xi(t)),
$$
where $M$ is the random variable defined in (\ref{M}).
\end{propos}
\medskip
We see that Theorem \ref{Kesten} ii) is a direct consequence of the
second estimate and of the tail estimate for $K(t)$. Theorem
\ref{doubleM}  is a consequence of the estimate i) and of the
estimate for $J$.
\medskip
\begin{proof}
{\it Step 1:} We first restrict the expectations to large values of
$H$. Let $h:\R_+\mapsto \R_+$  be any increasing function such that
\begin{eqnarray}\label{cond-ht1}
\lim_{t\to\infty} t^{-1}\ee^{h(t)}=0,
\\
\label{cond-ht2}
 h(t)\ge {9\over 10} \log t .
\end{eqnarray}
From Corollary \ref{ht}, we know that
\begin{eqnarray}\label{Hh}
\P^\Qi\left( \ee^H M_1 M_2\ge t, \; H\le
h(t)\right)=o(t^{-\kappa})=o(K(t)).
\end{eqnarray}
Hence, we can restrict ourselves to consider
$$
I_h(t):= \P^\Qi\left( \ee^H M_1 M_2\ge t\;|\; H\ge h(t)\right),
$$
$$
J_h(t):= \P^\Qi\left( \ee^H M_2 \ge t\;|\; H\ge h(t)\right),
$$
 \ali
{\it Step 2:} (Truncation of $M_1$, $M_2$). We need to truncate the
sums $M_1$ and $M_2$ so that they do not overlap. Under $\Qi(\cdot |
H\ge h(t))$ we consider the random variables
\begin{eqnarray}\label{tM1}
\tilde M_1:=\sum_{-\infty}^{t_1} \ee^{-V_k},
\end{eqnarray}
\begin{eqnarray}\label{tM2}
\tilde M_2:=\sum_{t_2}^\infty \ee^{V_k-S},
\end{eqnarray}
where
$$
t_1:=\inf\{k\ge 0,\; V_k\ge {1\over 3}\log t\}-1,
$$
$$
t_2:=\sup\{k\le T_H, \; V_k\le H-{1\over 3}\log t\}+1.
$$
Since $h(t)\ge {9\over 10}\log t$, we have
$$
0\le t_1<t_2\le T_H.
$$
Clearly, by the symmetry property of $\Qi$, $\tilde M_1$ and $\tilde
M_2$ have the same law under $\Qi(\cdot| H\ge h(t))$. (Observe that
the random variables $\tilde M_1$ and $\tilde M_2$ are implicitly
defined in terms of the variable $t$.)
\medskip
\begin{lem}\label{truncation}
Let $\xi$ be a positive real. There exists a constant $c_\xi>0$ such
that
$$
\P^\Qi\left( \tilde M_1\le \ee^{-\xi} M_1\;|\; H\ge h(t)\right)\le
\left\{\begin{array}{ll} c_\xi t^{{-\kappa/ 6}}&\hbox{ for
$\kappa\le 1,$}
\\
c_\xi t^{{-1/ 6}}&\hbox{ for $\kappa \ge 1$}.
\end{array}\right.
$$
\end{lem}
\medskip
\begin{proof}
We have, since $M_1\ge 1$
\begin{eqnarray*}
&& \P^\Qi\left( \tilde M_1\le \ee^{-\xi} M_1\;|\; H\ge h(t)\right)
\\
&\le & \P^\Qi\left( M_1-\tilde M_1\ge 1-\ee^{-\xi} \;|\; H\ge
h(t)\right)
\\
&\le &{1\over 1-\ee^{-\xi}} \E^\Qi\left( M_1-\tilde M_1 \;|\; H\ge
h(t)\right)
\\
&\le & c {\ee^{-\kappa h(t)}\over \P^\Qi(H\ge h(t))}\E^{\tQpos}
\Bigg( \sum_{k=t_1+1}^\infty \ee^{-Y_{k}}\Big( \sum_{{e_p\ge k,
\atop Y_{e_p}\ge h(t)}} \ee^{-\kappa (Y_{e_p}-h(t))}\Big)\Bigg),
\end{eqnarray*}
where in the last expression we used the result of Lemma
\ref{artificial2}, and the notation of the related section, and
where $c$ is a constant depending on $\xi$ and on the parameters of
the model. Using the fact that $\P^\Qi(H\ge h(t))\sim C \ee^{-\kappa
h(t)}$, when $t \to \infty,$ the Markov property and the fact that
$$
\E^{\tQpos}\Bigg(\sum_{{e_p\ge k,\atop \; Y_{e_p}\ge h(t)}}
\ee^{-\kappa (Y_{e_p}-h(t))}\Bigg)\le {1\over \P^{\tilde Q}(Y_n>0,
\forall n>0)( 1-\E^{\tilde Q}(\ee^{-\kappa Y_{e_1}}))},
$$
independently of $k$, we see that
\begin{eqnarray*}
\P^\Qi\Big( \tilde M_1\le \ee^{-\xi} M_1\;|\; H\ge h(t)\Big) &\le &
c \E^\tQpos \big( \sum_{k=t_1+1}^\infty \ee^{-Y_k}\big)
\\
&\le & c_\xi t^{-{\kappa\wedge 1 \over 6}},
\end{eqnarray*}
using the estimate of Lemma \ref{l:moment3}.
\end{proof}
\medskip
\noindent{\it Step 3:} (A small modification of the conditioning.)
We set
$$
\iii_h^{(t)}:=\iii\cap \{S\ge h(t)\}=\{V_k\ge 0\, , \, \forall k\le
0\}\cap\{S=H\}\cap \{S\ge h(t)\},
$$
the event by which we condition in $I_h(t)$, $J_h(t)$. We set
$$
\tIh:=\{S\ge h(t)\}\cap\{V_k\ge 0 \, , \,  \forall k\le 0\}\cap
\{V_k>0, \; \forall 0<k<T_{\third \log t}\},
$$
where
$$
T_{\third \log t}:=\inf\{ k\ge 0, \; V_k\ge \third\log t\}.
$$
Clearly, we have $ \Ih\subset \tIh $ and
$$
\P( \tIh\setminus \Ih\;|\; \tIh)\le c t^{-\kappa /3},
$$
for a constant $c>0$ depending only on the parameters of the model.
We set
$$
\tilde I_h(t):= \P^Q\left( \ee^H \tilde M_1 \tilde M_2\ge t\;|\;
\tIh \right),
$$
$$
\tilde J_h(t):= \P^Q\left( \ee^H \tilde M_2 \ge t\;|\; \tIh \right),
$$
$$
\tilde K_h(t):= \P^Q\left( \ee^H \ge t\;|\; \tIh \right).
$$
From Step 2 (Lemma \ref{truncation}) and Step 3, we see that we
have, for all $\xi>0,$ the following estimate
\begin{eqnarray}
\label{ItildeI} I_h(\ee^{2\xi}t)-c_\xi t^{-{\kappa\wedge 1\over
6}}\le \tilde I_h(t)\le I_h(t)+c t^{-{\kappa\over3}},
\\
\label{JtildeJ}
 J_h(\ee^{\xi}t)-c_\xi t^{-{\kappa\wedge 1\over  6}}\le \tilde
J_h(t)\le J_h(t)+c t^{-{\kappa\over3}}.
\end{eqnarray}
\ali {\it Step 4:} (The coupling strategy.)

 Let $(Y'_k)_{k\ge 0}$
and $(Y''_k)_{k\ge 0}$ be two independent processes with law
$$
\tilde Q(\cdot \;|\; Y_k>0 \, , \, 0<k \le T_{\third \log t}).
$$
Let us define, for all $u>0$, the hitting times
$$
T'_u:=\inf\{k\ge 0, \; Y'_k\ge u\}, \; \; T''_u:=\inf\{k\ge 0, \;
Y''_k\ge u\}.
$$
Set
$$
N'_0:=T'_{\third \log t}, \;\;\; N''_0:=T''_{\third\log t}.
$$
We couple the processes $(Y'_{N'_0+k})_{k\ge 0}$ and
$(Y''_{N_0''+k})_{k\ge 0}$ as in Durrett (cf \cite{durrett}, (4.3),
p. 204): we construct some random times $K'\ge N_0'$ and $K''\ge
N''_0$ such that
$$
\vert Y'_{K'}-Y''_{K''}\vert \le \xi,
$$
and such that $(Y'_{K'+k}-Y'_{K'})_{k\ge 0}$ and
$(Y''_{K''+k}-Y''_{K''})_{k\ge 0}$ are independent of the
$\sigma$-field generated by $Y'_0,\ldots ,Y'_{K'}$ and  $Y''_0,
\ldots , Y''_{K''}$. The method for this $\xi$-coupling is the
following: we consider some independent Bernoulli random variables
$(\eta'_i)_{i\in \N}$ and $(\eta''_i)_{i\in \N}$ (with
$\P(\eta'_i=1)=\P(\eta''_i=1)=\demi$) and we define
$$
(Z'_k)=(Y'_{N'_0+\sum_{i=1}^k \eta_i'}), \;\;\;\;
(Z''_k)=(Y''_{N''_0+\sum_{i=1}^k \eta_i''}).
$$
This extra randomization ensures that the process $(Z'_k-Z''_k)$ is
non arithmetic. Since its expectation is null, there exists a
positive random time for which  $Z'_k$ and $Z''_k$ are at a distance at most $\xi$
(cf the proof of Chung-Fuchs theorem (2.7), p. 188 and theorem
(2.1), p. 183 in \cite{durrett}). Then we define
$$
Y_k=\left\{\begin{array}{ll} Y'_k, \;&\hbox{when $k\le K',$}
\\
(Y''_{K''+(k-K')}-Y''_{K''})+Y'_{K'}, &\hbox{when $k> K'.$}
\end{array}\right.
$$
Clearly, by construction, since the processes $Y'$ and $Y''$ are no
longer conditioned when they reach the level $\third \log t$,
$(Y_k)_{k \ge 0}$ has the law
$$
\tilde Q(\cdot |\; Y_k>0, \; \forall \, 0<k< T_{\third\log t} ).
$$
We want that $Y'$ and $Y''$ to couple before they reach the level
$\demi \log t$, so we set
$$
\aaa=\{ K' <T'_{\demi \log t}\}\cap \{K''<T''_{\demi\log t}\}.
$$

Clearly, since the distribution of $Y'_{N'_0}-\third \log t$
converges (and the same for $Y''$, cf limit theorem (4.10), p. 370
in \cite{feller}) and since for all starting points $Y'_{N_0'}$ and
$Y''_{N''_0}$, $Z'$ and $Z''$ couple in a finite time almost surely,
we have the following result (whose proof is postponed to the end of
the section).

\medskip
\begin{lem}\label{couplage}
$$
\lim_{t\to \infty} \P (\aaa^c)=0.
$$
\end{lem}
\medskip

We set
$$
\eta(t):=\P(\aaa^c),
$$
and we choose $h(t)$ in terms of $\eta$ by
\begin{eqnarray}\label{choixh}
h(t)= (\log t+{1\over 2\kappa} \log \eta(t))  \vee ({9\over 10}\log
t) \vee ((1-{1\over 7\kappa})\log t),
\end{eqnarray}
where $\vee$ stands for the maximum of the three values. Clearly,
$h(t)$ satisfies the hypotheses (\ref{cond-ht1}), (\ref{cond-ht2}).

Consider now two independent processes $(W_k)_{k\ge 0}$ and
$(W'_k)_{k\ge 0}$ (and independent of $Y', Y''$) with the same law
$Q^{\le 0}$ (cf Subsection \ref{subsec:artificial}). Let $e$ be a strictly
increasing ladder time of $Y$ and define the process
$V(W,W',Y,e)=(V_k)_{k\in \Z}$ by
$$
\left\{ \begin{array}{l} (V_k)_{k\le 0}=(-W_{-k})_{k\le 0},
\\
\\
(V_k)_{k\ge 0}=(Y_0,\ldots ,Y_e, Y_e+W'_1, \ldots, Y_e+W'_k,\ldots).
\end{array}
\right.
$$
If $Y_e\ge h(t)$ then clearly $(V_k)_{k\in \Z}$ belongs to the event
$\tIh$, and the functional $ \tilde M_1 $ defined in (\ref{tM1})
depends only on $W$ and $Y'$; we denote it by $\tilde M_1(W,Y')$.
The functional $\tilde M_2$ depends only on $Y,W',e$; we denote it
by $\tilde M_2(Y,W',e)$. Using Lemma \ref{artificial}, we see that
\begin{eqnarray*}
\tilde I_h(t)={1\over \zzz_h(t)}\E\Big( \sum_{p=0}^\infty
\ee^{-\kappa Y_{e_p}} \indic_{Y_{e_p}\ge h(t)} \indic_{\tilde
M_1(W,Y')\tilde M_2(Y,W',e_p)\ee^{Y_{e_p}}\ge t}\Big),
\end{eqnarray*}
where $(e_p)_{p\ge 0}$ is the set of strictly increasing ladder times
of $Y$ (cf Subsection \ref{subsec:artificial}) and where $\zzz_h(t)$ is
the normalizing constant
$$
\zzz_h(t)= \E\Big( \sum_{p=0}^\infty \ee^{-\kappa Y_{e_p}}
\indic_{Y_{e_p}\ge h(t)}\Big).
$$
Clearly, $\zzz_h(t)\sim_{t\to\infty} c \ee^{-\kappa h(t)}$. The
variable $Y_{T_{h(t)}}-h(t)$ is indeed the residual waiting time of
the renewal process defined by the values of the process $Y$ at the
successive increasing ladder epochs. Hence, it converges in
distribution by the limit theorem (4.10) in (\cite{feller}, p. 370).

On the coupling event $\aaa$, we have
\begin{eqnarray*}
& Y''_{e_p-K'+K''}-\xi\le Y_{e_p}\le Y''_{e_p-K'+K''} +\xi,
\\
& \tilde M_2(Y,W',e_p)= \tilde M_2(Y'',W',e_p-K'+K''),
\end{eqnarray*}
for all ladder times $e_p$ such that $Y_{e_p}\ge h(t)$ (indeed
$h(t)\ge{9\over 10}\log t$) and where $\tilde
M_2(Y'',W',e_p-K'+K'')$ is the functional obtained from the
concatenation of the processes $Y''$ and $W'$ at time $e_p-K'+K''$,
as done for $\tilde M_2(Y,W',e_p)$. The first set of
inequalities implies that, on the coupling event $\aaa$, the set
$\{e_p-K'+K'', \; Y_{e_p}\ge h(t)\}$ is included in the set of
strictly increasing ladder times of $Y''$ larger than $h(t)-\xi$. So
we have
\begin{eqnarray*}
\tilde I_h(t)&\le &{\ee^{\kappa\xi}\over \zzz_h(t)}\E\Big(
\indic_\aaa\big(\sum_{p=0}^\infty \ee^{-\kappa Y''_{e''_p}}
\indic_{Y''_{e''_p}\ge h(t)-\xi} \indic_{\tilde M_1(W,Y')\tilde
M_2(Y'',W',e''_p)\exp(Y''_{e''_p})\ge t \ee^{- \xi}}\big)\Big)
\\
&& +{\ee^{-\kappa h(t)}\over \zzz_h(t)}\E\Big(
\indic_{\aaa^c}\big(\sum_{p=0}^\infty \ee^{-\kappa
(Y_{e_p}-h(t))}\indic_{Y_{e_p}\ge h(t)}\big)\Big),
\end{eqnarray*}
where $(e''_p)_{p \ge 0}$ denote the strictly increasing ladder times
for the process $Y''$.  Since the process $\{Y_{e_p},\; Y_{e_p}\ge
h(t)\}$ depends on the event $\aaa$ only through the value of
$Y_{T_{h(t)}}$, we see that the second term is less than or equal to
\begin{eqnarray}\label{Ac}
{1\over 1-\E^{\tilde Q}(\ee^{-\kappa Y_{e_1}})}{\ee^{-\kappa
h(t)}\over \zzz_h(t)} \P(\aaa^c)\le c\P(\aaa^c).
\end{eqnarray}
Now, the first term is lower than
\begin{eqnarray}\label{1term}
 \; \; \; \; \; \; \; \ee^{\kappa\xi}{\zzz_{h-\xi}(t)\over \zzz_h(t)}\P\big(\ee^{S''} \tilde
M_1(W,Y') \tilde M_2''\ge t \ee^{-\xi} \big)\le \ee^{3\kappa \xi}
\P\big(\ee^{S''} \tilde M_1(W,Y') \tilde M_2''\ge t \ee^{-\xi}\big),
\end{eqnarray}
for $t$ large enough (using the equivalent of $\zzz_h(t)$), where
$S''$ and $\tilde M_2''$ are relative to a process $V''$ independent
of $W,Y'$ and with law $Q(\cdot \;| \; \tilde \iii_{h-\xi}^{(t)})$.
Moreover, let us introduce $M''_2:=\sum_{k=0}^\infty
\ee^{V''_k-S''}.$ We need now to replace the truncated sum $\tilde
M_1$ by $M$. Using the fact that $\P(\exists k>0:
Y'_k\le 0)\le ct^{-\kappa/3}$, we see that
\begin{eqnarray}\label{1termbis}
\P\big(\ee^{S''} \tilde M_1(W,Y') \tilde M_2''\ge t \ee^{-\xi}
\big)& \le& \P\big( \ee^{S''} M_2'' M\ge t \ee^{-\xi}\big)
+ct^{-\kappa/6}
\\
&\le& \E\big( J_{h-\xi}(\ee^{-\xi}t/M)\big)+c't^{-\kappa/6},
\nonumber
\end{eqnarray}
the second inequality being a consequence of $\P( \tIh\setminus
\Ih\;|\; \tIh)\le c t^{-\kappa /3}$ and $M$ the random variable
defined in (\ref{M}) and independent of $V''$. Finally, considering
the choice made for $h(t)$ (cf (\ref{choixh})), we have
$$
t^{-{\kappa\wedge 1\over  6}} \P^\Qi (H\ge h(t)) =o(t^{-\kappa}),
$$
$$
\P(\aaa^c) \P^\Qi (H\ge h(t)) \le c t^{-\kappa}
\sqrt{\P(\aaa^c)}=o(t^{-\kappa}).
$$
Putting everything together (i.e., the estimates (\ref{Hh}),
(\ref{ItildeI}), (\ref{Ac}), (\ref{1term}), (\ref{1termbis}))
\begin{eqnarray*}
I(t)&\le & \P(H\ge h(t)) I_h(t)+o(t^{-\kappa})
\\
&\le&  \P(H\ge h(t))( \tilde I_h(\ee^{-2\xi} t)+c t^{-{\kappa\wedge
1\over  6}})+ o(t^{-\kappa})
\\
&\le & \P(H\ge h(t))(\ee^{3\kappa \xi}\E(J_{h-\xi}(\ee^{-3\xi}
t/M))+c\P(\aaa^c)) +o(t^{-\kappa})
\\
&\le & \ee^{3\kappa \xi} \P\left( R M\ge t \ee^{-3\xi}, \; H\ge
h(t)-\xi \right) +o(t^{-\kappa}),
\end{eqnarray*}
where $R$ and $M$ are independent processes with laws defined in
Section \ref{sec:main} (indeed, in the last inequality, $\P(H\ge
h(t))\P(\aaa^c) \le \sqrt{\P(\aaa^c)} t^{-\kappa}=o(t^{-\kappa})$).
Now, proceeding exactly as in Corollary \ref{ht}, we see that
$$
\P(RM\ge t, \; H<h(t)-\xi)=o(t^{-\kappa}),
$$
(indeed, the only difference is that $M_1$ is replaced by $M$ and
that $M$ and $R$ are independent). Finally, we proved that
$$
I(t)\le \ee^{3\kappa\xi} \E(J(\ee^{-3 \xi} t/M))+o(t^{-\kappa}).
$$
The lower estimate is similar. We first have, since the set
$\{e_p-K'+K'', \;\; Y_{e_p}\ge h(t)\} $ includes the set of strictly
increasing ladder times of $Y''$ larger than $h(t)+\xi$:
\begin{eqnarray*}
\tilde I_h(t)&\ge &{\ee^{-\kappa\xi}\over \zzz_h(t)}\E\bigg(
\indic_\aaa\big(\sum_{p=0}^\infty \ee^{-\kappa Y''_{e''_p}}
\indic_{Y''_{e''_p}\ge h(t)+\xi} \indic_{\tilde M_1(W,Y')\tilde
M_2(Y'',W',e''_p)\exp(Y''_{e''_p})\ge t \ee^{\xi}}\big)\bigg).
\end{eqnarray*}
Hence, by the same argument as above
\begin{eqnarray*}
\tilde I_h^{(t)}\ge \ee^{-3\kappa \xi} \P\big(\ee^{S''} \tilde
M_1(W,Y') \tilde M_2''\ge \ee^{\xi} t\big)+c\P(\aaa^c),
\end{eqnarray*}
where $S''$ and $\tilde M_2''$ are relative to a process $V''$
independent of $W$ and $Y'$ and with law $Q(\cdot\;| \; \tilde
I_{h+\xi}^{(t)})$.  Using, now the fact that $Y'_k>0$ for all $k>0$
with probability at least $1-ct^{-\kappa/3}$ and the fact that
$\tilde M_2\ge \ee^{-\xi} M_2$ with probability at least
$1-ct^{-\kappa/6}$, and the estimate on the tail of the sum $\sum
\ee^{-Y'_k}$ (of Subsection \ref{subsec:moment}) we see that
\begin{eqnarray*}
\tilde I_h^{(t)}\ge \ee^{-3\kappa \xi} \P\left(M \ee^{S''} M_2''\ge
\ee^{3\xi} t\right)+o(t^{-\kappa/6})+c\P(\aaa^c),
\end{eqnarray*}
where $M$ is the random variable defined in (\ref{M}) and
independent of $V''$. Then, we conclude as previously.

To prove the estimate on $J(t)$ and $K(t)$ we proceed exactly in the
same way: we first remark that by the property of time reversal (see
Lemma \ref{l:return}), we have
$$
J(t)=\P^\Qi (\ee^H M_1\ge t).
$$
The situation is then even simpler, we just have to decouple $M_1$
and $\ee^H$.
\end{proof}

\begin{proof} (of Lemma \ref{couplage}).
Denote by $F_{y',y''}(u)$ the probability that $Z'$ and $Z''$ couple
before the level $\third \log t+ u$ knowing that $Y'_{N_0'}=\third
\log t +y'$ and $Y''_{N_0''}=\third\log t +y''$. By the arguments
above, $F_{y',y''}(u)$ tends to 1 when $u$ tends to infinity. Let
$A>0;$ we first prove that this convergence is uniform in $y',y''$
on the compact set $y'\le A$, $y''\le A$. For this we consider the set
$\sss=(\N\cdot {\xi\over 4})\cap [0,A]$, and for $y'$, $y''$ in
$\sss\times \sss$ the function $\widehat F_{y',y''}(u),$ the
probability that $Z'$ and $Z'',$ starting from the points
$Y'_{N_0'}=\third \log t +y'$ and $Y''_{N_0''}=\third \log t +y'',$
couple at a distance $\xi/2$, before the level $\third \log t
+u-\xi$. Let
$$
\phi(u)=\inf_{y'\in \sss, \; y''\in \sss} \widehat F_{y',y''}(u).
$$
Clearly $\phi(u)\to 1$ when $u\to \infty$ and $F_{y',y''}(u)\ge
\phi(u)$, whenever $y'$ and $y''$ are in $[0,A]$. This implies that
$$
\liminf_{t\to \infty}\P(\aaa) \ge \liminf_{A\to \infty}\liminf_{t\to
\infty} \Big(\P(Y'_{N_0'}-\third \log t \le A)\Big)^2.
$$
Moreover, $\P^{\t Q}(Y'_k>0 \, , \, 0<k \le T_{\third \log t})\ge
\P^{\t Q}(Y'_k>0 \, , \, k \ge 0)>0$ implies
$$\P(Y'_{N_0'}-\third
\log t \ge A)=\P^{\t Q}(V_{T_{\third \log t}}-\third \log t \ge A \,
| \, V>0) \le c \P^{\t Q}(V_{T_{\third \log t}}-\third \log t \ge
A),
$$
where here $V$ is the canonical process under $\t Q.$ Therefore,
since $V_{T_{\third \log t}}-\third \log t $ converges in law (under
$\t Q$) to a finite random variable when $t$ tends to infinity (see
limit theorem (4.10), p. 370 in \cite{feller} or Example 4.4 part
II, page 214 in \cite{durrett}), this yields $\liminf_{t\to
\infty}\P(\aaa)=1.$
\end{proof}
\medskip

\section{Proof of Theorem \ref{Kesten},Theorem \ref{doubleM} and Theorem \ref{Kestengene}}
\label{sec:proof}
\begin{proof} (of Theorem \ref{Kesten}, $ii)$ and Theorem
\ref{doubleM}). Let $\xi>0$. By Proposition \ref{p:encadrement}, we
have, for all $A>0$ and for $t$ large enough,
\begin{eqnarray*}
J(t)\le \ee^{3\xi} \left(\E(K(\ee^{-2\xi} t M^{-1})\indic_{M\le
A})+\E(K(\ee^{-2\xi} t M^{-1})\indic_{M> A})\right).
\end{eqnarray*}
On the first term, for $t$ large enough, we can bound $K(\ee^{-2\xi} t M^{-1})$ from above
 by $({C_{I}\over \P^Q(H=S)}+\xi) (t
\ee^{-2\xi} M^{-1})^{-\kappa}$. For the second term we can use a
uniform bound $K(t)\le ct^{-\kappa}$. Thus we get
\begin{eqnarray*}
J(t)\le \ee^{3(1+\kappa)\xi}({C_{I} \over \P^Q(H=S)}+\xi)t^{-\kappa}
(\E( M^{\kappa})\indic_{M\le A})+c t^{-\kappa} \E(
M^{\kappa}\indic_{M> A})).
\end{eqnarray*}
Since $M^\kappa$ is integrable, letting $A$ tend to $\infty$, then
$\xi$ tend to 0, we get the upper bound
$$
\limsup_{t\to \infty} t^{\kappa} J(t)\le{C_{KI}\over \P^Q(H=S)}.
$$
For the lower bound it is the same. The proof of Theorem
\ref{doubleM} is the same: we use the estimate i) of Proposition
\ref{p:encadrement} and the tail estimate for $J$.
\end{proof}
\medskip
\begin{proof} (of Theorem \ref{Kesten}, $i)$).

Let us first recall (\ref{equivkesten}) and Theorem \ref{Kesten},
$ii),$ which tells us that
\begin{eqnarray}
\label{equivkestenigle+} Q(R>t \, ;\, H=S)=
\frac{C_{KI}}{t^{\kappa}}+o(t^{-\kappa}), \qquad t \to \infty,
\end{eqnarray}

\noindent where $C_{KI}=C_I \e(M^\kappa).$ Then, introducing
$$
KI:=\sum_{0 \le k \le T_{\r_-}} \ee^{V_k}, \qquad
O_1:=-V_{T_{\r_-}},$$ Theorem \ref{Kesten}, $i)$ is a consequence of
Theorem \ref{Kesten}, $ii)$ together with the two following lemmas.

\medskip
\begin{lem}
 \label{l:KItoKI} We have
\begin{eqnarray}
\label{equivkestenigle} Q(KI>t)=
\frac{C_{KI}}{t^{\kappa}}+o(t^{-\kappa}), \qquad t \to \infty.
\end{eqnarray}
\end{lem}
\medskip

\begin{proof} Firstly, observe that $KI \le R$ implies
$Q(KI>t \, ;\, H=S) \le Q(R>t \, ;\, H=S).$ Moreover, Corollary
\ref{ht} implies $Q(KI>t \, ;\, \ee^H=\ee^S\le
t^{2/3})=o(t^{-\kappa}),$ $t \to \infty,$ since $KI \le R.$
Furthermore, we have $0 \le Q(KI>t \, ;\, \ee^H>t^{2/3})-Q(KI>t \,
;\, \ee^H=\ee^S>t^{2/3}) \le Q(H \neq S \, ;\,
\ee^H>t^{2/3})=o(t^{-\kappa}),$ $t \to \infty.$ Therefore, we
obtain, when $t \to \infty,$
\begin{eqnarray}
\label{kitokieq1}
 Q(R>t \, ;\, H=S) \ge Q(KI>t \, ;\, \ee^H>t^{2/3})+
 o(t^{-\kappa}).
\end{eqnarray}

\noindent Since, by Corollary \ref{ht}, $Q(KI>t \, ; \, \ee^H\le
t^{2/3})=o(t^{-\kappa}),$ $t \to \infty,$ we get
\begin{eqnarray}
\label{kitokieq2} Q(KI>t \, ;\,
\ee^H>t^{2/3})=Q(KI>t)+o(t^{-\kappa}),
\end{eqnarray}

\noindent when $t \to \infty.$ Then, assembling (\ref{kitokieq1})
and (\ref{kitokieq2}) yields
\begin{eqnarray}
\label{kitokieq3} Q(R>t \, ;\, H=S)\ge Q(KI>t)+o(t^{-\kappa}),
\qquad t \to \infty.
\end{eqnarray}

On the other hand, observe that Corollary \ref{ht} implies that
$Q(R>t \, ;\, H=S)=Q(R>t \, ;\,
\ee^H=\ee^S>t^{2/3})+o(t^{-\kappa}),$ $t \to \infty.$ Moreover,
since we have $R=KI+\ee^{O_1}R',$ with $R'$ a random variable
independent of $KI$ and $O_1,$ having the same law as $R,$ we obtain
that $Q(R>t \, ;\, \ee^H=\ee^S>t^{2/3})\le Q_1 +Q_2,$ where
\begin{eqnarray*}
  Q_1 &:=& Q(KI \le t-t^{2/3} \, ; \, R'>t^{2/3} \, ;\, \ee^H>t^{2/3}),
  \\
  Q_2 &:=& Q(KI>t-t^{2/3} \, ; \, R >t \, ;\,
  \ee^H=\ee^S>t^{2/3}).
\end{eqnarray*}

\noindent Now, since $R'$ and $H$ are independent, we get $Q_1 \le
Q(\ee^H>t^{2/3}) Q(R'>t^{2/3})=o(t^{-\kappa}),$ $t \to \infty.$
Moreover, we easily have $Q_2 \le Q(KI>t-t^{2/3}).$ Therefore
\begin{eqnarray}
\label{kitokieq4} Q(R>t \, ;\, H=S)\le
Q(KI>t-t^{2/3})+o(t^{-\kappa}), \qquad t \to \infty.
\end{eqnarray}

\noindent Recalling (\ref{equivkestenigle+}) and assembling
(\ref{kitokieq3}) and (\ref{kitokieq4}) concludes the proof of Lemma
\ref{l:KItoKI}. \end{proof}

\medskip
\begin{lem}
 \label{l:KItoKesten}$C_{KI}$ satisfies
$$
C_{KI}=(1-\e^{Q}(\ee^{-\kappa O_1}))C_K.
$$
\end{lem}
\medskip

\begin{proof} First, observe that $Q(R>t)=Q(KI>t)+P_1+P_2,$
where
\begin{eqnarray*}
  P_1 &:=& Q(KI + \ee^{-O_1} R' >t \, ; \, t^{1/2}<KI \le t),
  \\
  P_2 &:=& Q(KI + \ee^{-O_1} R' >t \, ; \,  KI \le t^{1/2}),
\end{eqnarray*}

\noindent with $R'$ a random variable independent of  $KI$ and
$O_1,$ with the same law as $R.$

Now, let us prove that $P_1$ is negligible. Observe first that,
since $O_1\ge0$ by definition, we have $P_1 \le Q(R' >t-KI \, ; \,
t^{1/2}<KI \le t).$ Therefore $0 \le P_1 \le P_1'+P_1'',$ where
\begin{eqnarray*}
  P_1' &:=& Q(R' >t-KI \, ; \, t-t^{2/3} < KI \le t),
  \\
  P_1'' &:=& Q(R' >t-KI \, ; \, t^{1/2} < KI \le  t-t^{2/3}).
\end{eqnarray*}

\noindent Since $R'$ and $KI$ are independent, (\ref{equivkesten})
and (\ref{equivkestenigle}) yield $P_1'' \le Q(R' >t^{2/3}) Q(KI
>t^{1/2})=o(t^{\kappa}),$ $t \to \infty.$ Furthermore, we have
\begin{eqnarray*}
  P_1' &\le& Q(t-t^{2/3} < KI \le t)
  \\
       &\le& Q(KI >t-t^{2/3}) - Q(KI >t)
  \\
       &=& Q(KI >t) \left( \frac{Q(KI >t-t^{2/3})}{Q(KI >t)}-1  \right).
\end{eqnarray*}

\noindent Therefore (\ref{equivkestenigle}) implies
$P_1'=o(t^{-\kappa}),$ $t \to \infty.$ Then, we obtain
$P_1=o(t^{-\kappa}),$ $t \to \infty.$

Now, let us estimate $P_2.$ Observe that $\underline{P}_2\le P_2 \le
\overline{P}_2,$ where
\begin{eqnarray*}
  \underline{P}_2 &:=& Q(\ee^{-O_1} R' >t \, ; \, KI \le t^{1/2}),
  \\
  \overline{P}_2 &:=& Q(\ee^{-O_1} R' >t - t^{1/2}).
\end{eqnarray*}

\noindent Since $R'$ and $O_1$ are independent, (\ref{equivkesten})
yields
\begin{eqnarray}
\label{equivoverP2} \overline{P}_2=\frac{\e^Q(\ee^{-\kappa O_1})
C_{K}}{t^{\kappa}} + o(t^{-\kappa}), \qquad t \to \infty.
\end{eqnarray}

\noindent Therefore, it only remains to estimate $\underline{P}_2.$
Since $R'$ is independent of $KI$ and $O_1,$ we obtain for any
$\varepsilon>0$ and $t$ large enough,
$$
(1-\varepsilon) \, C_K \e^Q\bigg( {\bf 1}_{\{ KI \le t^{1/2} \}} \,
\frac{\ee^{-\kappa O_1}}{t^{\kappa}} \bigg) \le \underline{P}_2 \le
(1+\varepsilon) \, C_K \e^Q\bigg( {\bf 1}_{\{ KI \le t^{1/2} \}} \,
\frac{\ee^{-\kappa O_1}}{t^{\kappa}} \bigg).
$$

\noindent Moreover,
\begin{eqnarray*}
  \e^Q\left({\bf 1}_{\{ KI \le t^{1/2} \}} \,  \frac{\ee^{-\kappa O_1}}{t^{\kappa}} \right) &=&
  \frac{\e^Q(\ee^{-\kappa O_1})}{t^{\kappa}}-\e^Q\left( {\bf 1}_{\{ KI > t^{1/2} \}} \,
  \frac{\ee^{-\kappa O_1}}{t^{\kappa}}\right),
\end{eqnarray*}

\noindent and the second term on the right-hand side is less or
equal than $t^{-\kappa} Q(KI > t^{1/2})= o(t^{-\kappa}),$  $t \to
\infty.$ Thus
\begin{eqnarray}
\label{equivunderP2} \underline{P}_2=\frac{\e^Q(\ee^{-\kappa O_1})
C_{K}}{t^{\kappa}} + o(t^{-\kappa}), \qquad t \to \infty.
\end{eqnarray}

\noindent Assembling (\ref{equivoverP2}) and (\ref{equivunderP2})
yields ${P}_2=\frac{\e^Q(\ee^{-\kappa O_1}) C_{K}}{t^{\kappa}} +
o(t^{-\kappa}),$ $t \to \infty.$ Therefore, recalling
(\ref{equivkesten}), (\ref{equivkestenigle}) and
$Q(R>t)=Q(KI>t)+P_1+P_2,$ we obtain $C_{KI}=(1-\e^Q(\ee^{-\kappa
O_1}))C_K,$ which concludes the proof of Lemma \ref{l:KItoKesten}.
\end{proof}

\medskip

Since Theorem \ref{Kesten}, $ii)$ together with Lemma \ref{l:KItoKI}
and Lemma \ref{l:KItoKesten} yield $C_{KI}=C_I
\e^Q[M^\kappa]=(1-\e^{Q}[\ee^{-\kappa O_1}])C_K,$ we get $C_K=C_I
\e^Q[M^\kappa](1-\e^{Q}[\ee^{-\kappa O_1}])^{-1}.$ Now, recalling
that $C_I=(1-\e^{Q}[\ee^{-\kappa O_1}])C_F,$ this concludes the
proof of Theorem \ref{Kesten}, $i)$.
\end{proof}
\medskip
\begin{proof} (of Theorem \ref{Kestengene}).

The proof of Theorem \ref{Kestengene} is based on the same arguments
as in the proof of Theorem \ref{Kesten}, $i).$  We mainly have to
check analogous statements to Lemma \ref{prelim} and Corollary
\ref{ht}. Namely, we check that there exists $c>0$ such that
$$
\E^{Q^\iii}((M_1^B)^{\kappa+\frac{\varepsilon}{2}}\; |\; \lfloor
H\rfloor)\le c, \;\;\; Q^\iii \hbox{- a.s.,}
$$
where $M_1^B:=\sum_{k=-\infty}^{T_S}\ee^{-V_k}\tilde{B}_k.$ Using the
H\"older inequality instead of the Cauchy-Schwarz inequality in the
proof of Lemma \ref{prelim} we are led to check the integrability of
$(M^B)^{\kappa+\varepsilon}.$ This is used in the proof of
$$
\P^{Q^\iii}\left( R^B \ge t, \; H\le h(t)\right)=o(t^{-\kappa}),
$$
when $t$ tends to infinity, which is analogous to the proof of
Corollary \ref{ht} (in its $R$ version), choosing
$\eta=\kappa+\frac{\varepsilon}{2}.$

Now, it only remains to check the integrability of
$(M^B)^{\kappa+\varepsilon}.$ To this aim, we prove that
$\E^{\tQpos}((\sum_{k\ge 0}
\ee^{-V_k}\tilde{B}_k)^{\kappa+\varepsilon})<\infty,$ the case of
$\E^{Q^{>0}}((\sum_{k< 0}
\ee^{-V_k}\tilde{B}_k)^{\kappa+\varepsilon})$ being similar.

If $\kappa \ge 1,$ the Minkowski inequality yields
\begin{eqnarray}
\E^{\tQpos}\Big((\sum_{k\ge 0}
\ee^{-V_k}\tilde{B}_k)^{\kappa+\varepsilon}\Big) &\le& \left(
\sum_{k\ge 0} \E^{\tQpos}\Big((
\ee^{-V_k}\tilde{B}_k)^{\kappa+\varepsilon}\Big)^{\frac{1}{\kappa+\varepsilon}}
\right)^{\kappa+\varepsilon} \nonumber
\\
&\le& c \left( \sum_{k\ge 0} \E^{\tQpos}(
\ee^{-(\kappa+\varepsilon)V_k})^{\frac{1}{\kappa+\varepsilon}}
\right)^{\kappa+\varepsilon} \nonumber
\\
&\le& c \left( \sum_{k\ge 0} \E^{\tQpos}(
\ee^{-V_k})^{\frac{1}{\kappa+\varepsilon}}
\right)^{\kappa+\varepsilon}, \label{sum}
\end{eqnarray}
the second inequality being a consequence of the independence
between $(\tilde{B}_i)_{i \ge 0}$ and $(V_i)_{i \ge 0},$ while the
third inequality is due to the fact that $V_i \ge 0$ for $i \ge 0$
under $\tQpos$ together with $\kappa + \varepsilon\ge 1.$ Choosing
$p$ such that $p/(\kappa+\varepsilon)>1,$ let us write
$$
\E^{\tQpos}( \ee^{-V_k}) \le \frac{1}{k^p} +  \P^{\tQpos}(
\ee^{-V_k} \ge k^{-p}).
$$
Now, as in the proof of Lemma \ref{l:moment1}, since large
deviations do occur, we get from Cramer's theory, see
\cite{dembo-zeitouni98}, that the sequence $(\P^{\tQpos}( \ee^{-V_k}
\ge k^{-p}))_{k \ge 1}$ is exponentially decreasing. This yields
that the sum in (\ref{sum}) is finite.

If $\kappa<1,$ observe that we can restrict our attention to the
case where $\kappa+\varepsilon<1.$ Then, let us write
\begin{eqnarray*}
\E^{\tQpos}\Big((\sum_{k\ge 0}
\ee^{-V_k}\tilde{B}_k)^{\kappa+\varepsilon}\Big) &\le&
\E^{\tQpos}\Big(\sum_{k\ge 0}
(\ee^{-V_k}\tilde{B}_k)^{\kappa+\varepsilon}\Big)
\\
&\le& c \sum_{k\ge 0} \E^{\tQpos}( \ee^{-(\kappa+\varepsilon)V_k}),
\end{eqnarray*}
the second inequality being a consequence of the independence
between $(\tilde{B}_i)_{i \ge 0}$ and $(V_i)_{i \ge 0}.$ Now, the
conclusion is the same as in the case $\kappa \ge 1.$
\end{proof}

\section{Appendix}
\label{sec:appendix}

\subsection{Preliminaries' proofs} We give here the proofs of the claims from Section \ref{sec:preliminaries}.
\medskip

\noindent{\it Proof of Lemma \ref{l:moment1}.}
Using the Markov inequality, we get
$$\e^{\t Q}\Big(\sum_{k\geq0} \ee^{-V_k}\,|\, V\geq -L\Big) \le 1+ \sum_{k\geq1}{1\over k^2}+\sum_{k\geq1} \t
Q\Big(\ee^{-V_k}\geq {1\over k^2}\,|\, V\geq-L\Big) \ee^L.$$

Since $\P^{\t Q}(V \ge -L)\ge\P^{\t Q}(V \ge 0)>0,$ for all $L\ge0,$
$$\t
Q\Big(\ee^{-V_k}\geq {1\over k^2}\,|\, V\geq-L\Big)=\t Q(V_k\leq
2\log k\,|\, V\geq-L)\leq c \, \t Q(V_k\leq 2\log k).$$

Now, since large deviations do occur, we get, from Cramer's theory,
see \cite{dembo-zeitouni98}, that $ \e^{\t Q}(\log \rho_0)>0$
implies that the sequence $\t Q(V_k\leq 2\log k)$ is exponentially
decreasing.

The sum $\sum_{k\geq1} \t Q\big(\ee^{-V_k}\geq {1\over k^2}\,|\,
V\geq-L\big)$ is therefore bounded uniformly in $L$, and the result
follows.\qed
\bigskip

\noindent {\it Proof of Lemma \ref{l:moment2}.}
Let us treat, for more readability, the case
of the second moment. Observe first that
$$\e^{\t Q^{\ge0}}\bigg(\Big(\sum_{i\geq0} \ee^{-V_i}\Big)^2\bigg)
\leq 2\e^{\t Q^{\ge0}}\bigg(\sum_{i\geq0}\ee^{-2V_i}\Big(\sum_{j\geq
i} \ee^{-(V_j-V_i)}\Big)\bigg).$$ Applying the Markov property to
the process $V$ under $\t Q$ at time $i$, we get
$$\e^{\t Q^{\ge0}}\bigg(\Big(\sum_{i\geq0} \ee^{-V_i}\Big)^2\bigg)
\leq 2\e^{\t Q^{\ge0}}\bigg(\sum_{i\geq0} \ee^{-2V_i}\e^{\t
Q}\Big[\sum_{l\geq0}\ee^{-V'_l}\, |\, V' \geq -V_i\Big]\bigg),$$
where $V'$ is a copy of $V$ independent of $(V_k)_{0 \le k \le i}.$
Now, we use Lemma \ref{l:moment1} to get the upper bound
$$c \,  \e^{\t Q^{\ge0}}\Big(\sum_{i\geq0} \ee^{-2V_i}\times
\ee^{V_i}\Big) \leq c \,  \e^{\t Q^{\ge0}}\Big(\sum_{i\geq0}
\ee^{-V_i}\Big),$$ which is finite, again by applying Lemma
\ref{l:moment1}. This scheme is then easily extended  to higher
moments.
\qed
\bigskip

\noindent {\it Proof of Lemma \ref{l:moment3}.}
Let $\alpha\in[0,1]$ and define
$T_{(-\infty,-\alpha L]}:=\min\{i\geq0: V_i \leq -\alpha L\}$. Let
us write

$\begin{array}{rl}\ds\sum_{i\geq0}\ee^{-V_i}=&\ds\Big(\sum_{i\geq0}\ee^{-V_i}\Big) {\bf 1}_{\{V>-\alpha L\}}\\
&\ds+\bigg(\sum_{i=0}^{T_{(-\infty,-\alpha
L]}-1}\ee^{-V_i}+\sum_{i=T_{(-\infty,-\alpha L]}}^\infty
\ee^{-V_i}\bigg) {\bf 1}_{\{T_{(-\infty,-\alpha
L]}<\infty\}}.\end{array}$

Now, since $\t Q(V\geq -A)$ is uniformly bounded below, for
$A>0$, by $\t Q(V>0)>0,$ we obtain that $\e^{\t
Q}\big(\sum_{i\geq0}\ee^{-V_i}\, |\, V\geq -L\big)$ is less than or equal
to
\begin{eqnarray} &&c \,   \e^{\t Q}\big(\sum_{i\geq0}\ee^{-V_i}\,
|\, V\geq -\alpha L\big) \label{moment1}
\\
&+& c \,   \e^{\t Q}\bigg(\sum_{i<T_{(-\infty,-\alpha L]}}\ee^{-V_i} \,
; \, T_{(-\infty,-\alpha L]}<\infty \, ; \, V\geq -L\bigg) \nonumber
\\
&+&  c \,   \e^{\t Q}\bigg(\sum_{i\geq T_{(-\infty,-\alpha
L]}}\ee^{-V_i} \, ;\, T_{(-\infty,-\alpha L]}<\infty\, ; \, V\geq
-L\bigg). \nonumber
\end{eqnarray}

\noindent Lemma \ref{l:moment1} bounds the first term in
(\ref{moment1}) from above by $c  \ee^{\alpha L},$ for all $L>0.$ Furthermore,
$i<T_{(-\infty,-\alpha L]}$ implies $\ee^{-V_i} \le \ee^{\alpha L}.$
Therefore, $c  \ee^{\alpha L} \e^{\t Q}(T_{(-\infty,-\alpha L]} {\bf
1}_{\{T_{(-\infty,-\alpha L]}<\infty\}})$ is an upper bound for the
second term in (\ref{moment1}), which is treated as follows,
\begin{eqnarray*} \e^{\t Q}\big(T_{(-\infty,-\alpha L]}{\bf 1}_{\{T_{(-\infty,-\alpha
L]}<\infty\}}\big)&\leq& \sum_{k\geq 0} k\t
Q\big(T_{(-\infty,-\alpha L]}=k\big)
\\
&\leq& \sum_{k\geq 0} k \t Q\big(V_k\leq-\alpha L\big)
\\
&\leq& \sum_{k\geq0}k \ee^{-k \theta \t I(-{\alpha L\over k})}
\ee^{-k (1-\theta) \t I(-{\alpha L\over k})},
\end{eqnarray*}
where $0<\theta<1$ and $\t I$ denotes the rate function associated
with $\t P$ which is positive convex and admits a unique minimum  on
$\R_+.$ We can therefore bound below all the terms  $\t
I(-{\alpha h\over k})$ by $\t I(0)>0.$ Moreover, a more
sophisticated result yields $\sup_{x \le 0} \t I(x)/x \le - \kappa$
(see definition of $\kappa$ and formula $(2.2.10)$ in
(\cite{dembo-zeitouni98}, p.~$28$)). Therefore, we obtain
$$ \e^{\t Q}\big(T_{(-\infty,-\alpha L]}{\bf 1}_{\{T_{(-\infty,-\alpha
L]}<\infty\}}\big)\leq  \ee^{- \theta \kappa \alpha L}
\sum_{k\geq0}k \ee^{-k (1-\theta) \t I(0)} \le c \,  \ee^{- \theta
\kappa \alpha L}.$$ As a result, the second term in (\ref{moment1})
is bounded by $c  \ee^{(1-\theta \kappa \alpha)L},$ for all $L>0.$

Finally, concerning the third term in (\ref{moment1}), we have that
\begin{eqnarray*}
&&c \,  \e^{\t Q}\bigg(\sum_{i\geq T_{(-\infty,-\alpha
L]}}\ee^{-V_i}\, ;\, T_{(-\infty,-\alpha L]}<\infty \, ; \, V\geq
-L\bigg)
\\
&\le& c \,  \e^{\t Q}\bigg(\ee^{-V_{T_{(-\infty,-\alpha L]}}}
\sum_{i\geq T_{(-\infty,-\alpha L]}}
\ee^{-(V_i-V_{T_{(-\infty,-\alpha L]}})}\, ; \, T_{(-\infty,-\alpha
L]}<\infty \,; \, V\geq -L\bigg)
\\
&\le& c \, \e^{\t Q}\bigg(\ee^{-V_{T_{(-\infty,-\alpha L]}}} {\bf
1}_{\{T_{(-\infty,-\alpha L]}<\infty\}} \e^{\t Q}\Big( \sum_{i\geq
0} \ee^{-V'_i} \,| \, V' \ge -(L+ V_{T_{(-\infty,-\alpha L]}})
\Big)\bigg),
\end{eqnarray*}

\noindent where $V'_i:=V_{T_{(-\infty,-\alpha
L]}+i}-V_{T_{(-\infty,-\alpha L]}},$ for $i \ge 0.$ The last
inequality is a consequence of the strong Markov property applied at
$T_{(-\infty,-\alpha L]},$ which implies that $(V'_i, \, i \ge 0)$
is a copy of $(V_i, \, i \ge 0)$ independent of $(V_i,\, 0 \le i \le
T_{(-\infty,-\alpha L]}).$ Then, Lemma \ref{l:moment1} yields that
the third term in (\ref{moment1}) is less than
\begin{eqnarray*}
&& c \,  \e^{\t Q}\bigg(\ee^{-V_{T_{(-\infty,-\alpha L]}}} {\bf
1}_{\{T_{(-\infty,-\alpha L]}<\infty\}} \ee^{L+
V_{T_{(-\infty,-\alpha L]}}} \bigg)
\\
&\le& c \,  \ee^L \t Q(T_{(-\infty,-\alpha L]}<\infty) \le c \,
\ee^{(1-\kappa \alpha)L}.
\end{eqnarray*}
Since $\theta <1$ implies $1-\theta \kappa \alpha \ge 1- \kappa
\alpha,$ we optimize the value of $\alpha$ by taking $\alpha=-\alpha
\kappa \theta+1,$ i.e. $\alpha=1/(1+\kappa \theta)$. As a result, we
get already a finer result than Lemma \ref{l:moment1} with a bound
$\ee^{L\over 1+\kappa \theta}$ instead of $\ee^L$.

Now, the strategy is to use this improved estimation instead of
Lemma \ref{l:moment1} and repeat the same procedure.  In that way, we
obtain recursively a sequence of bounds, which we denote by $c  \ee^{u_n
L}$. The first term in (\ref{moment1}) is bounded by $c \ee^{\alpha
u_n L}$ whereas the second and the third term are still bounded
respectively by $c \ee^{(1-\kappa \alpha \theta) L}$ and $c
\ee^{(1-\kappa \alpha) L}.$

Optimizing in $\alpha$ again, one chooses $\alpha u_n=-\alpha \kappa
\theta+1,$ i.e. $\alpha={1\over u_n+\kappa \theta}.$ The new
exponent is therefore $u_{n+1}=\alpha u_n={u_n\over u_n+\kappa
\theta}.$ Thus, the sequence $u_n$ is monotone and converges to a
limit satisfying $l={l\over l+\kappa \theta}$ . For $\kappa
\theta\leq1,$ the limit is $l=1-\kappa \theta$ and for $\kappa
\theta\geq1,$ the limit is 0. Since this result holds for any
$0<\theta<1,$ it concludes the proof of Lemma \ref{l:moment3}.
\qed
\bigskip

\noindent {\it Proof of Lemma \ref{l:return}.}
Let $\phi$ be a positive test function. We have
\begin{eqnarray*}
&&\E^\Qi(\phi((V_{T_H}-V_{T_H-k})_{k\ge 0}))
\\
&=& \sum_{p=0}^\infty
\E^\Qi(\indic_{T_H=p}\phi((V_{p}-V_{p-k})_{k\ge 0}))
\\
&=& {1\over \P^Q(\iii)} \sum_{p=0}^\infty \E^Q(\indic_{\{V_k\ge 0\,
, \, \forall k\le 0\}}\indic_{\{V_k\le V_p \, , \, \forall k\ge
p\}}\indic_{\{0<V_k<V_p \, , \, \forall 0<k<p\}}
\phi((V_{p}-V_{p-k})_{k\ge 0})).
\end{eqnarray*}
By construction we have, for all $p\ge 0,$
$$
(V_p-V_{p-k})_{k\in \Z}\equilaw (V_k)_{k\in \Z}.
$$
This implies that
\begin{eqnarray*}
&&\E^\Qi(\phi((V_{T_H}-V_{T_H-k})_{k\ge 0}))
\\
&=&{1\over \P^Q(\iii)} \sum_{p=0}^\infty \E^Q(\indic_{\{V_k\ge 0\, ,
\, \forall k\le 0\}}\indic_{\{V_k\le V_p \, , \,\forall k\ge
p\}}\indic_{\{0<V_k<V_p \, , \, \forall 0<k<p\}} \phi((V_k)_{k\ge
0}))
\\
&=& \E^\Qi(\phi((V_k)_{k\in \Z})). \end{eqnarray*}
\qed
\bigskip

\noindent {\it Proof of Lemma \ref{artificial}.}
Let $\psi, \theta$ be positive test functions. Let us compute
\begin{eqnarray*}
&&\E^Q\left(\psi( (V_{T_S+k}-V_{T_S})_{k\ge 0})\theta ((V_0, \ldots,
V_{T_S})) \right)
\\
&=& \sum_{p=0}^\infty \E^Q\left( \indic_{T_S=p}\psi(
(V_{p+k}-V_{p})_{k\ge 0})\theta ((V_0, \ldots,V_{p}))\right)
\\
&=& \sum_{p=0}^\infty \E^Q\left(\indic_{\{V_k<V_p \, , \, \forall
0\le k<p\}}\indic_{\{V_k\le V_p \, , \, \forall k\ge p\}}\psi(
(V_{p+k}-V_{p})_{k\ge 0})\theta ((V_0, \ldots, V_{p}))\right)
\\
&=& \sum_{p=0}^\infty \E^Q\left( \indic_{\{V_k<V_p\, , \, \forall
k<p\}}\theta ((V_0, \ldots, V_{p}))\right) \E^Q\left(
\indic_{\{V_k\le 0\, , \,\forall k\ge 0\}}\psi( (V_{k})_{k\ge
0})\right),
\end{eqnarray*}
using the Markov property at time $p$. The second term is equal to
$$
\P^Q(V_k\le 0\, , \, \forall k\ge 0) \, \E^\Qneg \left(\psi(
(V_{k})_{k\ge 0})\right).
$$
Let us now consider only the first term. Using the Girsanov property
of $Q$ and $\tilde Q$ we get
\begin{eqnarray*}
 \sum_{p=0}^\infty \E^Q\left( \indic_{\{V_k<V_p\, , \, \forall
k<p\}}\theta ((V_0, \ldots, V_{p}))\right) &=& \sum_{p=0}^\infty
\E^{\tilde Q}\left( \indic_{\{V_k<V_p\, , \, \forall
k<p\}}\ee^{-\kappa V_p} \theta ((V_0, \ldots, V_{p}))\right)
\\
&=& \sum_{p=0}^\infty \E^{\tilde Q}\left( \ee^{-\kappa V_{e_p}}
\theta ((V_0, \ldots, V_{e_p}))\right),
\end{eqnarray*}
where $(e_p)_{p \ge 0}$ are the strictly increasing ladder times of
$(V_k \, , \, k \ge 0)$ as defined in Subsection \ref{subsec:artificial}. The last formula is
exactly the one we need, and also implies that
$$
{1\over \zzz} =\P^Q(V_k\le0 \, , \, \forall k\ge 0)=1-\E^{\tilde
Q}(\ee^{-\kappa V_{e_1}}),
$$
(which can also be obtained directly).
\qed
\bigskip

\noindent {\it Proof of Lemma \ref{artificial2}.}
Let $\Psi$ be a positive test function. Thanks to the previous lemma,
we have
\begin{eqnarray*}
&&\E^\Qi(\Psi(V_0, \ldots ,V_{T_H}))
\\
&=& {1\over \P^Q(H=S)}\E^Q(\indic_{H=S}\Psi(V_0, \ldots ,V_{T_H}))
\\
&=& {1\over \zzz \P^Q(H=S)}\sum_{p=0}^\infty \E^{\tilde Q}
(\indic_{Y_k>0\, , \, \forall 0<k\le e_p } \ee^{-\kappa Y_{e_p}}
\Psi(Y_0, \ldots ,Y_{e_p}))
\\
&=& {1\over \zzz \P^Q(H=S)}\sum_{p=0}^\infty \E^{\tilde  Q}
(\indic_{Y_k>0\, , \, \forall k>0 } {1\over \P^{\tilde
Q}(V_k>-Y_{e_p}\, , \, \forall k>0)} \ee^{-\kappa Y_{e_p}} \Psi(Y_0,
\ldots ,Y_{e_p}))
\\
&\le & {1\over \zzz \P^Q(H=S) \P^{\tilde Q}(V_k>0\, , \,\forall
k>0)}\sum_{p=0}^\infty \E^{\tilde  Q} (\indic_{Y_k>0\, , \, \forall
k>0 } \ee^{-\kappa Y_{e_p}} \Psi(Y_0, \ldots ,Y_{e_p}))
\\
&\le & {1\over  \P^Q(H=S) \P^{\tilde Q}(V_k>0\, , \,\forall
k>0)}\E^{\hat Q^{>0}} (\Psi(Y_0, \ldots ,Y_{\Theta})),
\end{eqnarray*}
using the Markov property at time $e_p$ in the fourth line. This is
exactly what we want.
\qed
\bigskip

\subsection{A Tauberian result}
\label{subsec:tauberian}
\medskip
\begin{corrol}\label{Tauberian}
Let $h:\R_+\rightarrow \R_+$ be such that
$$
\lim_{\lambda \to 0} \lambda \ee^{h(\lambda)}=0, \;\;\;
\lim_{\lambda\to 0} h(\lambda)=\infty.
$$
Then, for $\kappa<1,$
$$
\E^Q \Big( 1-{1\over 1+\lambda Z}\; |\; \iii_{h}^{(\lambda)} \Big)
{\sim}  {1 \over \P^Q(H\ge h(\lambda))}\; {\pi\kappa\over \sin(\pi
\kappa)} C_U \lambda^\kappa,
$$
when $\lambda\to 0$, where $\iii_{h}^{(\lambda)}$ is the event
$$
\iii_{h}^{(\lambda)}=\iii\cap\{H\ge h(\lambda)\}=\{V_k\ge 0\, , \,
\forall k\le 0\}\cap \{H=S\ge h(\lambda)\}.
$$
\end{corrol}
\medskip

\begin{proof}
Clearly, we have
$$
\E^Q \Big( 1-{1\over 1+\lambda Z}\; |\; \iii_{h}^{(\lambda)} \Big) =
{\P^Q(H=S)\over \P^Q(H=S\ge h(\lambda))} \E^\Qi \Big( \indic_{H\ge
h(\lambda)}\big(1-{1\over 1+\lambda Z}\big)\Big).
$$
Since $\P^Q(H=S\ge h(\lambda))\sim \P^Q(H\ge h(\lambda))$ we
consider now
$$
\E^\Qi \Big( \indic_{H\ge h(\lambda)}\big(1-{1\over 1+\lambda
Z}\big)\Big)
$$
 We will omit in the following the reference to the law
$\Qi$, and simply write $\E$ for the expectation with respect to
$\Qi$. We have
\begin{eqnarray}\label{un+Z}
&&\E\Big( \indic_{H\ge h(\lambda)} \big(1-{1\over 1+\lambda
Z}\big)\Big)
\\
&=&\E\Big( \indic_{Z\ge  \ee^{h(\lambda)}} \big(1-{1\over 1+\lambda
Z}\big)\Big)-\E\Big( \indic_{\ee^H < \ee^{h(\lambda)}\le Z}
\big(1-{1\over 1+\lambda Z}\big)\Big). \nonumber
\end{eqnarray}
For $\kappa<1,$ the second term can be bounded by
\begin{eqnarray*}
\E\Big( \indic_{\ee^H < \ee^{h(\lambda)}\le Z} \big(1-{1\over
1+\lambda Z}\big)\Big) &\le & \sum_{p=0}^{\lfloor h(\lambda)\rfloor}
\E\Big( \indic_{\intH=p} {\lambda Z\over 1+\lambda Z}\Big)
\\
&=&
 \sum_{p=0}^{\lfloor h(\lambda)\rfloor} \E\Big(
\indic_{\intH=p}\, \E\Big( {\lambda Z\over 1+\lambda Z}\; |\;
\intH=p\Big)\Big)
\\
&\le & \sum_{p=0}^{\lfloor h(\lambda)\rfloor}
 \E\Big(
\indic_{\intH=p}{c \lambda \ee^p \over 1+ c \lambda \ee^p}\Big),
\end{eqnarray*}
where, in the last inequality, we used the Jensen inequality and
Corollary \ref{prelim}, and where $c$ denotes a constant independent
of $\lambda$ (which may change from line to line). Now, since
$\P(\intH=p)\le c \ee^{-\kappa p}$ for a positive constant $c$, we
get that
\begin{eqnarray*}
\E\Big( \indic_{\ee^H < \ee^{h(\lambda)}\le Z} \big(1-{1\over
1+\lambda Z}\big)\Big) &\le&  c \lambda \sum_{p=0}^{\lfloor
h(\lambda)\rfloor} \ee^{(1-\kappa) p} \le  c' \lambda
\ee^{(1-\kappa))h(\lambda)} \\
&\le& c \lambda^\kappa (\lambda
\ee^{h(\lambda)})^{1-\kappa}=o(\lambda^{\kappa}),
\end{eqnarray*}
for $\kappa<1,$ since $\lambda \ee^{h(\lambda)} \to 0,$ $\lambda \to
0.$

By integration by parts, we see that the first term of (\ref{un+Z})
is equal to
\begin{eqnarray*}
&&\E\Big( \indic_{Z\ge h(\lambda)} \big(1-{1\over 1+\lambda
Z}\big)\Big)
\\
&=&\Big[ {\lambda z\over 1+\lambda z} \P(Z\ge z)
\Big]_{\ee^{h(\lambda)}}^\infty +\int_{\ee^{h(\lambda)}}^\infty
{\lambda \over (1+\lambda z)^2} \P(Z\ge z) \d z.
\end{eqnarray*}
The first term is lower than
$$
c \lambda \ee^{(1-\kappa) h(\lambda)}=c\lambda^\kappa (\lambda
\ee^{h(\lambda)})^{1-\kappa}= o(\lambda^\kappa),
$$
for $\kappa<1$. For the second term, let us suppose first that
$$
h(\lambda)\rightarrow \infty.
$$
We can estimate $\P(Z\ge z)$ by
$$
({C_U\over \P^Q(H=S)}-\eta)z^{-\kappa} \le \P(Z\ge z)\le ({C_U\over
\P^Q(H=S)}+\eta) z^{-\kappa},
$$
for any $\eta,$ when $\lambda$ is sufficiently small. Hence we are
led to compute the integral
$$
\int_{\ee^{h(\lambda)}}^\infty {\lambda\over 1+\lambda z}
z^{-\kappa} \d z =\lambda^\kappa \int_{\lambda \ee^{h(\lambda)}\over
1+\lambda \ee^{h(\lambda)}}^1 x^{-\kappa}(1-x)^\kappa \d x,
$$
(making the change of variables $x=\lambda z/(1+\lambda z)$). For
$\kappa<1$ this integral converges, when $\lambda\to 0,$ to
$$
\Gamma(\kappa+1)\Gamma(-\kappa+1)={\pi \kappa\over \sin(\pi
\kappa)}.
$$
\end{proof}
\noindent\Rm: Let us make a final remark useful for
\cite{enriquez-sabot-zindy}. If we truncate the series $M_1$ on the
right and on the left when $V_k$ reaches the level $A>0$, and if we
truncate $M_2$ when $H-V_k$ reaches the level $A$ then the results
of Theorem \ref{doubleM} and Corollary \ref{Tauberian} remain valid
just by replacing in the tail estimate $M$ by the the random walk $M$ truncated
at level $A$. More precisely, let $A>0$ and consider
$$
\overline M_1 =\sum_{k=t_1^-}^{t_1^+} \ee^{-V_k}, \;\;\; \overline
M_2 =\sum_{k=t_2^-}^{t_2^+} \ee^{V_k-H}
$$
where
$$
t_1^-=\sup \{k \le 0 ,\; V_k\ge A\}, \;\;\; t_1^+=\inf \{k\ge 0 ,\;
V_k\ge A\}\wedge T_H
$$
$$
t_2^-=\sup \{k \le T_H ,\; H-V_k\ge A\}\vee 0, \;\;\;
 t_2^+=\inf
\{k\ge T_H ,\; V_k\ge A\}
$$
then the results of Theorem \ref{doubleM} and Corollary
\ref{Tauberian} remain valid when we consider $\overline Z= \ee^H
\overline M_1 \overline M_2$ instead of $Z$, if we replace in the
tail estimate $M$ by $\overline M=\sum_{t_-}^{t_+} \ee^{-V_k}$ where
$t_-$ and $t_+$ are the hitting times of the level $A$ on the left
and on the right. Indeed, in the proof of Theorem \ref{doubleM} we
see that considering the truncated $\overline M_1$ and $\overline
M_2$ only simplifies the proof: we don't need to truncate $M_1$ and
$M_2$ as we did. In particular, it implies that in Corollary
\ref{Tauberian} we can truncate $M_1$ and $M_2$ at a level
$\underline h(\lambda)\le h(\lambda)$: if $\underline h(\lambda)$
tends to $\infty,$ we have exactly the same result.

\bigskip
\bigskip
\noindent {\bf Acknowledgements}
The authors would like to thank two anonymous referees for  many useful comments
and   suggestions.

\bigskip

\end{document}